\documentclass[preprint,12pt]{elsarticle}
\usepackage{}

\usepackage{amsmath}
\usepackage{cases}
\usepackage{color}
\usepackage{hyperref}
\usepackage{amssymb}
\usepackage{graphicx}
\usepackage{graphics}
\usepackage{subfigure}
\usepackage{appendix}
\usepackage{textcomp}
\usepackage{bm}
\usepackage{url}

\makeatletter
  \newcommand\figcaption{\def\@captype{figure}\caption}
  \newcommand\tabcaption{\def\@captype{table}\caption}
\makeatother

\biboptions{numbers,sort&compress}

\begin{document}

\begin{center}

{\LARGE On the homotopy analysis method for backward/forward-backward stochastic differential equations}

\vspace{0.3cm}

Xiaoxu Zhong $^3$, Shijun Liao $^{1,2,3}$  \footnote{Corresponding author.  Email address: sjliao@sjtu.edu.cn}

\vspace{0.3cm}

$^1$ State Key Laboratory of Ocean Engineering, Shanghai 200240, China\\
$^2$ Ministry-of-Education Key Laboratory of Scientific and Engineering Computing\\
$^3$ School of Naval Architecture, Ocean and Civil Engineering\\
Shanghai Jiao Tong University, Shanghai 200240, China

 \end{center}

\hspace{-0.6cm}{\bf Abstract}  {
\em In this paper, an analytic approximation method for highly nonlinear equations, namely the homotopy analysis method (HAM), is employed to solve some backward stochastic differential equations (BSDEs) and forward-backward stochastic differential equations (FBSDEs),  including one with high dimensionality (up to 12 dimensions).  By means of the HAM, convergent series solutions can be quickly obtained with high accuracy for a FBSDE in a 6 dimensional case, within less than $1\%$ CPU time used by a currently reported numerical method for the same case [34].  Especially, as  dimensionality  enlarges, the increase of computational complexity for the HAM  is not as dramatic as this numerical method.   All of these demonstrate the validity and high efficiency  of the HAM for the backward/forward-backward stochastic differential equations  in science, engineering and finance.
}

\vspace{0.3cm}

\hspace{-0.6cm}{\bf Key Words} forward-backward stochastic differential equations,  \sep BSDEs, \sep  FBSDEs, \sep analytic approximation, \sep  homotopy analysis method

\renewcommand{\theequation}{1.\arabic{equation}}
\setcounter{equation}{0}
\section{Introduction}
The backward stochastic differential equations (BSDEs) are widely existed in scientific and engineering fields, such as stochastic control, stock markets, chemical reactions, and so on. The general form of BSDEs read
\begin{eqnarray}
\left\{
\begin{split}
-\mathrm{d}y_t&=&f(t, y_t, z_t)\mathrm{d}t-z_t\;\mathrm{d}W_t, \;\;\;\;0\leq t\leq T,\\
y_T&=&\xi,~~~~~~~~~~~~~~~~~~~~~~~~~~~~~~~~~~~~~~~~~~~ \label{1:general:bsde}
\end{split}
\right.
\end{eqnarray}
in which $W_t=(W_t^1,\cdot\cdot\cdot, W_t^d)^*$ is a $d$-dimensional Brownian motion defined in some complete, filtered probability space $(\Omega, {\cal F}, P, \{ {\cal F}_t\}_{0\leq t \leq T})$, $T>0$ denotes the deterministic terminal time, and the notation $^*$ is the transpose operator for matrix or vector. In finance, $y_t$ and $z_t$ denote the dynamic value of the replicating and hedging portfolios, respectively.

In 1990, the existence and uniqueness of solutions of nonlinear BSDEs were proved by Pardoux and Peng \cite{Pardoux1990}, which laid the theoretical foundation of the forward-backward stochastic differential equations (FBSDEs).  In 1991, Peng \cite{Peng19911, Peng19912} further found the relationship between FBSDEs and a kind of second-order quasi-linear parabolic partial differential equation. Since then, extensive researches of the BSDEs and FBSDEs have been done \cite{Peng1994, Pardoux1999, Peng1992, Peng1993, Peng1999, Peng1997The}. Besides, many interesting properties and applications of BSDEs were presented \cite{Peng1997}. So, how to efficiently obtain accurate solutions of BSDEs or FBSDEs becomes a significant problem.

In the complete, filtered probability space $(\Omega, {\cal F}, P, \{ {\cal F}_t\}_{0\leq t \leq T})$, the general form of FBSDEs read
\begin{eqnarray}
\left\{
\begin{split}
\mathrm{d}x_t&=&b(t, x_t, y_t, z_t)\mathrm{d}t+\sigma(t, x_t, y_t, z_t)\mathrm{d}W_t,~~\\
-\mathrm{d}y_t&=&f(t, x_t, y_t, z_t)\mathrm{d}t-z_t\mathrm{d}W_t, \;\;\;\;0\leq t\leq T,\\
x_0&=&\epsilon,\;\;\;\;\;\;\;\;y_T=\xi,~~~~~~~~~~~~~~~~~~~~~~~~~~~~~ \label{1:general:fbsde}
\end{split}
\right.
\end{eqnarray}
where $\epsilon\in {\cal F}_0$, $\xi\in {\cal F}_T$, $W_t$ is a $d$-dimensional Brownian motion; $b: \Omega\times[0, T]\times \mathbb{R}^d \times \mathbb{R}^m \times \mathbb{R}^{m\times r}\rightarrow \mathbb{R}^d$; $\sigma: \Omega\times[0, T]\times \mathbb{R}^d \times \mathbb{R}^m \times \mathbb{R}^{m\times r}\rightarrow \mathbb{R}^{d\times r}$; $f: \Omega\times[0, T]\times \mathbb{R}^d \times \mathbb{R}^m \times \mathbb{R}^{m\times r}\rightarrow \mathbb{R}^{m}$; $x_t \in \mathbb{R}^d$;  $y_t \in \mathbb{R}^m$ and $z_t \in \mathbb{R}^{m\times r}$ are unknown stochastic processes needed to be solved.

From the general Feynman-Kac formula \cite{Peng19912, Peng1997}, and also based on some smooth assumptions \cite{Pardoux1999}, the solution $(y_t, z_t)$ of (\ref{1:general:fbsde}) can be represented as
\begin{equation}
y_t=u(t, x_t),\;\;\;\;z_t=\sigma(t,x_t) \nabla_{x}u(t,x_t),  \label{1:fbsde:solution:relation}
\end{equation}
in which $u(t,x)$ is the classical solution of the partial differential equation
\begin{equation}
\frac{\partial u}{\partial t}+\sum_{i=1}^{d}b_i\frac{\partial u}{\partial x_i}+\frac{1}{2}\sum_{i, j=1}^{d}[\sigma \sigma^*]_{ij}\cdot\frac{\partial^2 u}{\partial x_i \partial x_j}+f(t, x, u, \sigma \nabla_x u)=0, \label{1:fbsde:corr:diff}
\end{equation}
subjects to the boundary condition
\begin{equation}
u(T,x)=\xi.  \label{1:general:fbsde:boundary}
\end{equation}

In finance, there is an another important type of stochastic differential equations, namely the second-order forward-backward stochastic differential equations (2FBSDEs). As an extension of the FBSDEs, the 2FBSDEs were first introduced and studied by Cheridito, Soner, Touzi and Victoir \cite{Cheridito2007Second}, and the relations between the 2FBSDEs and fully non-linear parabolic PDEs were revealed \cite{Cheridito2007Second, Soner2012Wellposedness}. Since then, variety of studies regarding the 2FBSDEs are carried out \cite{Tao2015High}. In the complete filtered probability space $(\Omega, {\cal F}, P, \{ {\cal F}_t\}_{0\leq t \leq T})$, the general form of the coupled 2FBSDEs read
\begin{eqnarray}
\left\{
\begin{split}
\mathrm{d}x_t&=&b(t, \Theta_t)\mathrm{d}t+\sigma(t, \Theta_t)\mathrm{d}W_t,~~~~~~~~~\\
-\mathrm{d}y_t&=&f(t, \Theta_t)\mathrm{d}t-z_t\mathrm{d}W_t, ~~~~~~~~~~~~~~~~\\
\mathrm{d}z_t&=&A_t\mathrm{d}t+\Gamma_t\mathrm{d}W_t,\;\;\;\;0\leq t\leq T,~~~~~\\
x_0&=&\epsilon,\;\;\;\;\;\;\;\;y_T=\xi,~~~~~~~~~~~~~~~~~~~~~~~ \label{1:general:2fbsde}
\end{split}
\right.
\end{eqnarray}
 in which $\Theta_t = (x_t, y_t, z_t, A_t, \Gamma_t)\in \mathbb{R}^m \times \mathbb{R} \times \mathbb{R}^d \times {\cal S}^d$ is an unknown stochastic process needed to be solved; ($\Omega$, ${\cal F}$, $P$) is a given probability space; $W_t$ is a $d$-dimensional Brownian motion; ${\cal S}^d$ is the set of all $d\times d$ real-valued symmetric matrices; $b: \Omega\times[0, T]\times \mathbb{R}^m \times \mathbb{R} \times \mathbb{R}^{d}\times {\cal S}^d\rightarrow \mathbb{R}^m$; $\sigma: \Omega\times[0, T]\times \mathbb{R}^m \times \mathbb{R} \times \mathbb{R}^{d}\times {\cal S}^d\rightarrow \mathbb{R}^{m\times d}$; $f: \Omega\times[0, T]\times \mathbb{R}^m \times \mathbb{R} \times \mathbb{R}^{d}\times {\cal S}^d\rightarrow \mathbb{R}$.

 Under some assumptions \cite{Cheridito2007Second, Soner2012Wellposedness}, the solution ($y_t$, $z_t$, $A_t$, $\Gamma_t$) of (\ref{1:general:2fbsde}) can be represented as
 \begin{eqnarray}
 \left\{
 \begin{split}
 y_t&=u(t,x_t),\;\;\;\;\;\;\;\;z_t= \sigma(t,x_t)  \nabla_{x}u(t,x_t), \\
 \Gamma_t &= \left(\sigma  \nabla_x(\sigma \nabla_x u)\right)(t,x_t), \;\;\;A_t=\left({\cal \tilde{L}}(\sigma \nabla_x u)\right)(t,x_t), \label{1:general:2fbsdes:boundary}
 \end{split}
 \right.
 \end{eqnarray}
 where ${\cal \tilde{L}}$ is a second-order elliptic differential operator
 \begin{equation}
 {\cal \tilde{L}}\phi=\frac{\partial \phi}{\partial t}+b(t,x)\nabla_x\phi+\frac{1}{2}\mathrm{tr}\left(\sigma(t,x)\sigma^*(t,x)\nabla_x^2\phi\right) \label{ellip:ope}
 \end{equation}
 with $\nabla_x\phi=(\partial_{x_1}\phi,\cdots, \partial_{x_m}\phi)$, and $\nabla_x^2\phi$ being the Hessian matrix \cite{Binmore2002Calculus} of $\phi$ with respect to the spatial variable $x$.

In the last twenty years, many numerical methods are proposed to solve the BSDEs, FBSDEs and 2FBSDEs. In 1994, Ma, Protter and Yong \cite{Ma1994} proposed a four step scheme to numerically solve the corresponding parabolic partial differential equations of FBSDEs. Based on the four step scheme, some numerical algorithms \cite{Douglas1996Numerical, Ma1999, Ma2002, Milstein2007} are proposed. In 2006, a time-space discretization scheme was proposed by Delarue and Menozzi \cite{Delarue2006A} for quasi-linear parabolic PDEs, which weakens the regularity assumptions required in \cite{Douglas1996Numerical}. In 2008, Delarue and Menozzi \cite{Fran2008An} improved this scheme by introducing a interpolation procedure. Apart from solving the parabolic partial differential equations, some numerical schemes are proposed to directly solve the BSDEs and FBSEDs \cite{Bally1997Approximation, Chevance1997Numerical, Bouchard2004Discrete, Zhang2004A, Gobet2005A, Bender2007A, Zhao2009Error, Zhao2014}. In 1999, Peng \cite{Peng2002A} proposed a linear approximation algorithm to solve the Chow's lagrangean one-dimensional BSDEs \cite{Chow1997Dynamic}. In 2006, Zhao, Chen and Peng \cite{Zhao2006A} proposed a $\theta$-scheme for BSDEs with high accuracy, and they gave a rigorous \emph{proof} that this scheme is of first-order convergence for general $\theta$, and especially is of second-order when $\theta=\frac{1}{2}$. Unfortunately, most of the above-mentioned numerical methods are not very efficient for high-dimensional FBSDEs (i.e. dimension exceeds 3), since the computational complexity increases dramatically as the dimension increases. To the best of our knowledge, up to now, quite few numerical schemes are proposed for high-dimensional FBSDEs \cite{Gobet2005A, Guo2012A, Fu2016Efficient}.

In this paper, the BSDEs, FBSDEs, 2FBSDEs and high-dimensional FBSDEs (i.e. up to 12-dimension) are successfully solved by means of an analytic approximation technique for highly nonlinear differential equations, namely the homotopy analysis method (HAM) \cite{liaoPhd, Liaobook, liaobook2, KV2012}. Unlike perturbation technique, the HAM is independent of any small/large physical parameters.  Besides, the HAM provides us great freedom to choose equation-type and solution expression of the high-order approximation equations. Especially, there is a convergence-control parameter $c_{0}$ in the series solutions, which provides us a convenient way to guarantee the convergence of series solutions gained by the HAM. It is these merits that distinguish the HAM from other analytic approaches, and thus enable the HAM to be successfully applied to many complicated problems with high nonlinearity \cite{Abbasbandy2006, Liang2010, Ghotbi2011, Nassar2011, Aureli2014, Duarte2015, Nagarajaiah2015, xu2012JFM, Liu2014JFM}. Note that the HAM was successfully applied to give an analytic approximation with much longer expiry for the optimal exercise boundary of an American put option than perturbation approximations\cite{Cheng2010An}. In this paper, by choosing proper auxiliary linear operator and convergence-control parameter $c_{0}$, convergent results of the BSDEs, FBSDEs, 2FBSDEs and high dimensional FBSDEs (up to 12-dimension) are efficiently obtained by means of the HAM. Our HAM approximations agree well with the exact solutions. All of these demonstrate the validity and potential of the HAM for the backward/forward-backward stochastic differential equations.

\renewcommand{\theequation}{2.\arabic{equation}}

\setcounter{equation}{0}
\section{The HAM for BSDEs}

In this section, we employ the HAM to solve the backward stochastic differential equations (BSDEs). According to (\ref{1:general:fbsde})-(\ref{1:general:fbsde:boundary}), if we let $x_t$ be a standard Brownian motion (i.e. $b=0$, $\sigma=1$), then the corresponding PDE related to the BSDEs reads
\begin{equation}
\frac{\partial u}{\partial t}+\frac{1}{2}\sum_{i=1}^{d}\frac{\partial^2 u}{\partial x_i^2}+f(t, u, \nabla_x u)=0, \;\;\;\;0 \leq t \leq T, \label{2:bsde:corr:diff}
\end{equation}
subjects to the boundary condition
\begin{equation}
u(T,x)=\xi,  \label{2:general:bsde:boundary}
\end{equation}
in which $d$ denotes the dimensionality. The solutions of the BSDEs satisfy
\begin{equation}
y_t=u(t, W_t),\;\;\;\;z_t=\nabla_{x}u(t,W_t). \label{2:bsde:solution:relation}
\end{equation}
\subsection{An example of one-dimensional BSDE}
At first, we consider the following BSDE with both $y_t$ and $W_t$ in one dimension
\begin{eqnarray}
\left\{
\begin{split}
-\mathrm{d}y_t&=&(-y_{t}^{3}+2.5y_{t}^{2}-1.5y_t)\mathrm{d}t-z_t\mathrm{d}W_t, \;\;\;\;0\leq t\leq T,\\
y_T&=&\frac{\mathrm{exp}(W_T+T)}{\mathrm{exp}(W_T+T)+1},~~~~~~~~~~~~~~~~~~~~~~~~~~~~~~~~~~~~ \label{21:case1:bsde}
\end{split}
\right.
\end{eqnarray}
with the exact solutions:
\begin{equation}
y_t=\frac{\mathrm{exp}(W_t+t)}{\mathrm{exp}(W_t+t)+1},\quad\quad\quad
z_t=\frac{\mathrm{exp}(W_t+t)}{(\mathrm{exp}(W_t+t)+1)^2}. \label{21:case1:solution}
\end{equation}
Without loss of generality, we choose $T=1$. Since $W_0=0$, the exact solutions at the initial moment are: $y_0=1/2$, $z_0=1/4$.
According to (\ref{2:bsde:corr:diff})-(\ref{2:bsde:solution:relation}), the corresponding PDE reads
\begin{equation}
\frac{\partial u}{\partial t}+\frac{1}{2}\frac{\partial ^2 u}{\partial x^2}-u^3+\frac{5}{2}u^2-\frac{3}{2}u=0,\;\;\;\;0\leq t\leq 1,  \label{21:pde:gov}
\end{equation}
subjects to the boundary condition
\begin{equation}
u(1,x)=\frac{\mathrm{e}^{x+1}}{\mathrm{e}^{x+1}+1}. \label{21:pde:boundary}
\end{equation}
Set
\begin{equation}
\theta=\frac{\mathrm{e}^{x+1}}{\mathrm{e}^{x+1}+1}.  \label{21:transform}
\end{equation}
Then $\theta\in (0,1)$, and we have
\begin{equation}
y_t=u(t,\theta),\;\;\;\;\;\;\;\;z_t=\theta(1-\theta)\frac{\partial u(t,\theta)}{\partial \theta}.
\end{equation}
Eqs.~(\ref{21:pde:gov}) and (\ref{21:pde:boundary}) become
\begin{eqnarray}
{\cal N}[u] &=& \frac{\partial u}{\partial t}+\frac{\theta}{2} (1-\theta)(1-2\theta)\frac{\partial u}{\partial \theta}+\frac{1}{2}\theta^2(1-\theta)^2\frac{\partial^2 u}{\partial \theta^2}-u^3 \nonumber \\
&+&\frac{5}{2}u^2-\frac{3}{2}u=0, \label{21:tran:pde:gov}
\end{eqnarray}
subjects to
\begin{equation}
u(t,\theta)\Big{|}_{t=1}=\theta, \label{21:tran:pde:boundary}
\end{equation}
with the exact solution
\begin{equation}
u(t,\theta)=\frac{\theta\mathrm{e}^t}{\theta \mathrm{e}^t+(1-\theta)\mathrm{e}}.
\end{equation}
Here, $\cal N$ denotes a nonlinear operator.

Since $t\in[0,1]$ is a finite interval, it is natural to express $u(t,\theta)$ in a power series
\begin{equation}
u(t,\theta)=\sum_{m=0}^{+\infty}a_m(\theta)\cdot t^m, \label{21:solution:exp}
\end{equation}
where $a_m(\theta)$ is coefficient to be determined. This provides us the so-called ``solution expression" of $u(t,\theta)$ in the frame of the HAM.

Let $\varphi_0(t,\theta)$ denote the initial guess of $u(t,\theta)$, which satisfies the boundary condition (\ref{21:tran:pde:boundary}). Moreover, let ${\cal L}$ denote an auxiliary linear operator with property ${\cal L}[0]=0$, $c_{0}$ a non-zero auxiliary parameter, called the convergence-control parameter, and $q\in [0,1]$ the embedding parameter for constructing a homotopy, respectively. Then we construct a family of differential equations
\begin{eqnarray}
\begin{split}
(1-q){\cal L}\bigg{[}\Phi(t,\theta,q)-\varphi_{0}(t,\theta)\bigg{]}&=c_{0} \; q  \;  {\cal N} \bigg{[} \Phi(t,\theta,q)\bigg{]}, \label{21:zeroth:equ}
\end{split}
\end{eqnarray}
subject to the boundary condition
\begin{equation}
\Phi(t,\theta,q)\Big{|}_{t=1}=\theta,  \label{21:zeroth:boundary}
\end{equation}
where the nonlinear operator $\cal N$ is defined by (\ref{21:tran:pde:gov}).  Note that  $\Phi(t,\theta,q)$ corresponds to the unknown $u(t,\theta)$, as mentioned below.

When $q=0$, due to the property ${\cal L}[0]=0$, Eqs.~(\ref{21:zeroth:equ}) and (\ref{21:zeroth:boundary}) have the solution
\begin{equation}
\Phi(t,\theta,0)=\varphi_0(t,\theta). \label{21:q0}
\end{equation}
When $q=1$, Eqs.~(\ref{21:zeroth:equ}) and (\ref{21:zeroth:boundary}) are equivalent to the original equations (\ref{21:tran:pde:gov}) and (\ref{21:tran:pde:boundary}), provided
\begin{equation}
\Phi(t,\theta,1)=u(t,\theta).  \label{21:q1}
\end{equation}
So, as the embedding parameter $q$ increases from $0$ to $1$, $\Phi(t,\theta,q)$ varies (or deform) continuously from $\varphi_0(t,\theta)$ to $u(t,\theta)$. Eqs.~(\ref{21:zeroth:equ}) and (\ref{21:zeroth:boundary}) are called the zeroth-order deformation equations.

Using (\ref{21:q0}), we have the Maclaurin power series with respect to the embedding parameter $q$:
\begin{equation}
\Phi(t,\theta,q)=\sum_{m=0}^{+\infty}\varphi_m(t,\theta) \cdot q^m,  \label{21:power:series}
\end{equation}
where
\begin{equation}
\varphi_m(t,\theta)={\cal D}_{m}\Big{[}\Phi(t,\theta,q)\Big{]},
\end{equation}
in which
\begin{equation}
{\cal D}_{m}[f]=\frac{1}{m!}\frac{\partial^{m}f}{\partial q^{m}}\bigg{|}_{q=0} \label{def:D}
\end{equation}
is called the $m$th-order homotopy-derivative \cite{Liaobook, liaobook2} of $f$.

In general, the convergence radius of a power series is finite. Fortunately, in the frame of the HAM, we have great freedom to choose the auxiliary linear operator $\cal L$ and the so-called convergence-control parameter $c_{0}$. Assume that all of them are so properly chosen that the power series (\ref{21:power:series}) is convergent at $q=1$. Then, according to (\ref{21:q1}), we have the so-called homotopy-series solution
\begin{equation}
u(t,\theta)=\sum_{m=0}^{+\infty}\varphi_{m}(t,\theta).
\end{equation}

Substituting (\ref{21:power:series}) into the zeroth-order deformation equations (\ref{21:zeroth:equ}) and (\ref{21:zeroth:boundary}), and then equating the like-power of $q$, we have the so-called $m$th-order deformation equations
\begin{equation}
{\cal L}\left[\varphi_m-\chi_m\varphi_{m-1}\right] = c_{0} \; \delta_{m-1}(t,\theta), \label{21:high:eqs}
\end{equation}
subject to the boundary condition
\begin{equation}
\varphi_{m}(t,\theta)\Big{|}_{t=1}=0, \label{21:high:boundary}
\end{equation}
where
\begin{eqnarray}
\delta_{n}(t,\theta) &=& {\cal D}_{n} \bigg{\{} {\cal N}[\Phi(t,\theta,q)] \bigg{\}} \nonumber\\
&=&\frac{\partial \varphi_{n}}{\partial t}+\frac{1}{2}\theta(1-\theta)(1-2\theta)\frac{\partial \varphi_{n}}{\partial \theta}+\frac{1}{2}\theta^2(1-\theta)^2\frac{\partial^2 \varphi_{n}}{\partial \theta^2} \nonumber \\
&& -\sum_{i=0}^{n}\sum_{j=0}^{n-i}\varphi_i \varphi_j \varphi_{n-i-j}+\frac{5}{2}\sum_{i=0}^{n}\varphi_i \varphi_{n-i}-\frac{3}{2}\varphi_{n},
\end{eqnarray}
and
\begin{equation}
    \chi_{m}=\left\{
    \begin{array}{ll}
    0,  \;\;m\leq1,\\
    1,  \;\;m>1. \end{array}\right.\label{def:chi}
\end{equation}
Here, ${\cal D}_n$ is defined by (\ref{def:D}).   To satisfy the condition (\ref{21:tran:pde:boundary}), we choose
\begin{equation}
\varphi_{0}(t,\theta)=\theta
\end{equation}
as the initial guess of $u(t,\theta)$.  Since there are no boundary conditions for $\theta$,  we simply choose such an auxiliary linear operator
\begin{equation}
{\cal L}[\varphi]=\frac{\partial \varphi}{\partial t}. \label{linear:operator}
\end{equation}
Note that, using the auxiliary linear operator (\ref{linear:operator}), the solution $\varphi_{m}(t,\theta)$ of the {\em linear} high-order deformation equations (\ref{21:high:eqs}) and (\ref{21:high:boundary})  can be easily obtained step by step, say,
\begin{equation}
\varphi_m(t,\theta) = \chi_m \; \varphi_{m-1}(t,\theta) + c_0 \; \int_1^t \delta_{m-1}(z,\theta) \text{d} z, \label{solution:integration}
\end{equation}
starting from $m=1$, especially by means of computer algebra software such as Mathematica.  The $M$th-order homotopy-approximation of $u(t,\theta)$ reads
\begin{equation}
\tilde{\varphi}_M(t,\theta)=\sum_{m=0}^{M}\varphi_{m}(t,\theta),
\end{equation}
which gives us the initial values:
\begin{equation}
\tilde{y}_0=\tilde{\varphi}_M(0,\theta)\Big{|}_{\theta=\frac{\mathrm{e}}{\mathrm{e}+1}},\;\;\;\;\;\;\tilde{z}_0=\Big{[} \theta(1-\theta)\frac{\partial \tilde{\varphi}_M(0,\theta)}{\partial \theta} \Big{]}\bigg{|}_{\theta=\frac{\mathrm{e}}{\mathrm{e}+1}}. \nonumber
\end{equation}
In other words,  in the frame of the HAM,  the nonlinear partial differential equation (\ref{21:tran:pde:gov}) with the boundary condition (\ref{21:tran:pde:boundary}) can be solved just by means of integrations with respect to the time $t$.   

Since the exact solution is known, we define the following squared residual error to characterize the global error between the homotopy approximation and the exact solution:
\begin{equation}
\tilde{{\cal E}}=\int_{0}^{1}\int_{0}^{1}\Big{[}\tilde{\varphi}_M(t,\theta)-u(t,\theta)\Big{]}^2\mathrm{d}\theta\mathrm{d}t.
\end{equation}
Besides,  we can always define such a squared residual error
\begin{equation}
{\cal E}=\int_{0}^{1}\int_{0}^{1}\Big{\{}{\cal N}[\tilde{\varphi}_M(t,\theta)]\Big{\}}^2\mathrm{d}\theta\mathrm{d}t,
\end{equation}
where ${\cal N}$ is defined by (\ref{21:tran:pde:gov}).
Obviously, the smaller the ${\cal E}$ or $\tilde{{\cal E}}$, the more accurate the HAM approximation.

\begin{figure}[t]
    \begin{center}
        \begin{tabular}{cc}
            \includegraphics[width=3in]{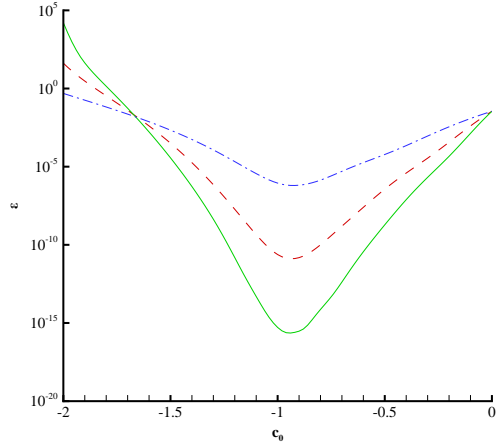}
        \end{tabular}
    \caption{The squared residual error versus $c_0$ of Eq.~(\ref{21:zeroth:equ}). Dash-dotted line: 5th-order approximation; Dashed line: 10th-order approximation; Solid line: 15th-order approximation.} \label{case1:c0}
    \end{center}
\end{figure}

The curves of ${\cal E}$ versus $c_0$ is as shown in Fig.~\ref{case1:c0}, which indicates that for any $c_0\in[-1.4,-0.4]$, ${\cal E}$ decreases as the order increases, i.e. the homotopy-series is convergent within the region $c_0\in[-1.4,-0.4]$, but the optimal $c_0$ (i.e. minimum ${\cal E}$) is near to $-1$. So, for the sake of simplicity, we choose $c_0=-1$, and the homotopy approximations are
\begin{equation}
\left\{
\begin{split}
\tilde{\varphi}_1(t,\theta)&=t\theta+\theta^2-t\theta^2;\\
\tilde{\varphi}_2(t,\theta)&=\frac{1}{2}\theta+\frac{1}{2}t^2\theta-\frac{1}{2}\theta^2+2t\theta^2
-\frac{3}{2}t^2\theta^2+\theta^3-2t\theta^3+t^2\theta^3;\\
\tilde{\varphi}_3(t,\theta)&=\frac{1}{3}\theta+\frac{1}{2}t\theta+\frac{1}{6}t^3\theta+
\frac{2}{3}\theta^2-\frac{3}{2}t\theta^2+2t^2\theta^2-\frac{7}{6}t^3\theta^2\\
&\quad-\theta^3+4t\theta^3-5t^2\theta^3+2t^3\theta^3+\theta^4-3t\theta^4+3t^2\theta^4-t^3\theta^4;\\
\cdots\quad&
\end{split}
\right.
\end{equation}

\begin{figure}[!t]
    \begin{center}
        \begin{tabular}{cc}
            \includegraphics[width=3in]{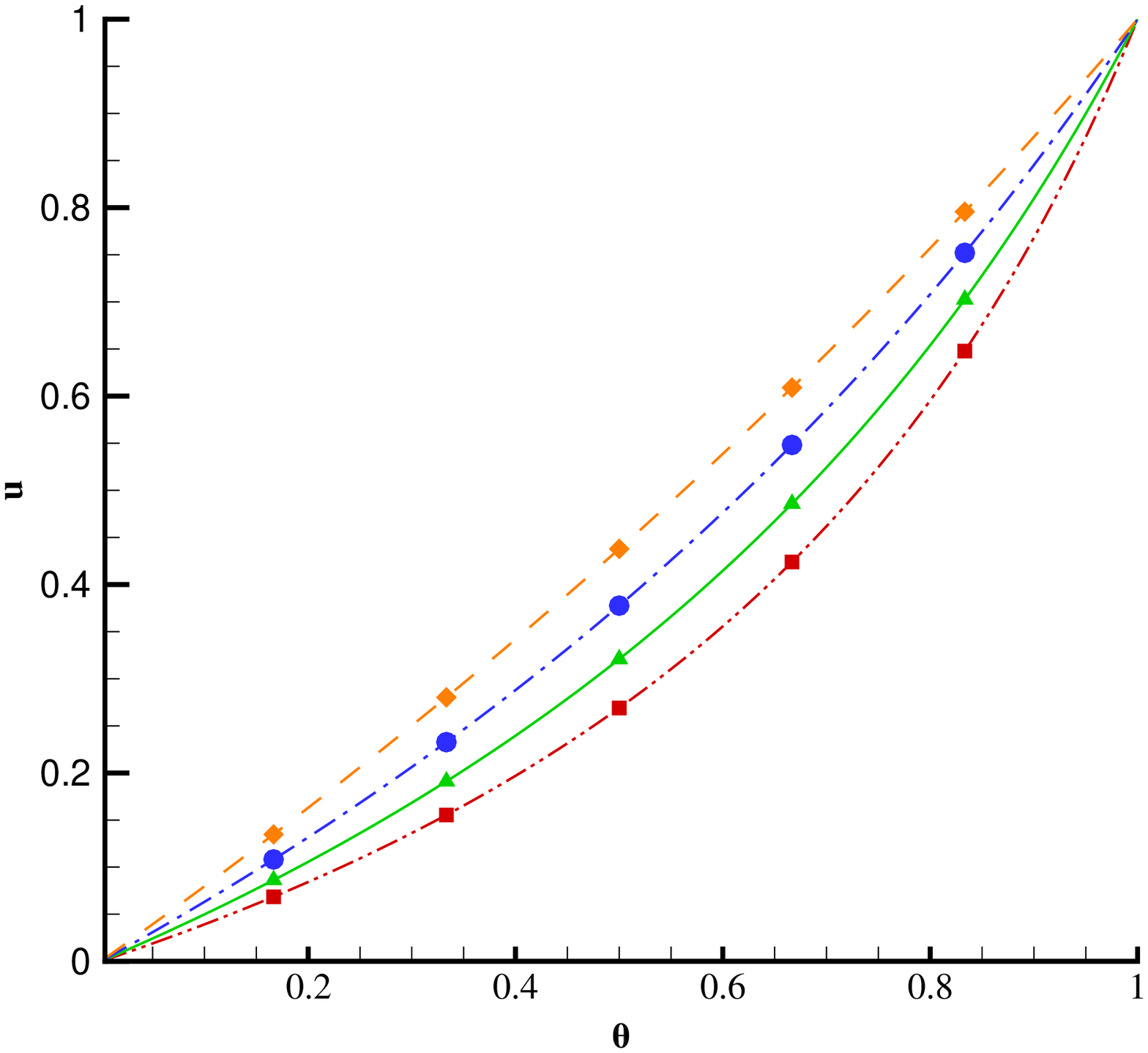}
        \end{tabular}
    \caption{Comparison between the homotopy approximations and the exact solution $u$ of Eqs.~(\ref{21:tran:pde:gov}) and (\ref{21:tran:pde:boundary}) when $t=0, 0.25, 0.5, 0.75$. Dash-double-dotted line: homotopy approximation at $t=0$; Solid line: homotopy approximation at $t=0.25$; Dashed-dotted line: homotopy approximation at $t=0.5$; Dashed line: homotopy approximation at $t=0.75$. Square: exact solution at $t=0$; Triangle up: exact solution at $t=0.25$; Circle: exact solution at $t=0.5$; Diamond: exact solution at $t=0.75$.} \label{case1:compare}
    \end{center}
\end{figure}

Note that rather accurate approximations at the level of $10^{-9}$ can be obtained in just 2.2 seconds CPU times\footnote{All examples considered in this paper are computed using a laptop with a core i7 3.60GHz process and 8GB memory.}, as shown in Table \ref{case1}. Besides, as shown in Fig.~\ref{case1:compare}, our homotopy approximations agree well with the exact solutions. All of these indicate that the HAM is valid and rather efficient for the BSDEs.

\begin{table}[!htb]
\tabcolsep 0pt
\caption{The squared residual error, the used CPU time (second) and the relative error between the homotopy approximations $\tilde{y}_0$, $\tilde{z}_0$ and the exact solutions $y_0$, $z_0$ of Eq.~(\ref{21:case1:bsde}), respectively. Here, $m$ denotes the order of approximation.}
\vspace*{-12pt}\label{case1}
\begin{center}
\def\temptablewidth{1\textwidth}
{\rule{\temptablewidth}{1pt}}
\begin{tabular*}{\temptablewidth}{@{\extracolsep{\fill}}ccccc}
$m$  &$\tilde{{\cal E}}$ &   $\tilde{y}_0-y_{0}$ &$\tilde{z}_0-z_{0}$    & time (s)   \\
\hline
3  & $1\times10^{-6}$ &  $-5\times10^{-3}$ & $-8\times10^{-4}$    &0.1     \\
6  & $8\times10^{-10}$ &  $8\times10^{-5}$ & $2\times10^{-4}$    &0.2       \\
9 &   $5\times10^{-13}$ & $3\times10^{-7}$& $-1\times10^{-5}$      &0.4        \\
12 &  $4\times10^{-16}$ &  $-8\times10^{-8}$& $5\times10^{-7}$    &1.0          \\
15 & $6\times10^{-19}$ &  $3\times10^{-9}$& $-8\times10^{-9}$     &2.2          \\
\end{tabular*}
{\rule{\temptablewidth}{1pt}}
 \end{center}
 \end{table}

 \subsection{An example of two-dimensional BSDE} %
 
In this case, we consider the following BSDE of $y_t$ in two dimensions
\begin{equation}
\left\{
\begin{split}
&-\left(
\begin{matrix} \mathrm{d}y_t^{(1)} \\ \mathrm{d}y_t^{(2)} \end{matrix}
\right)=
\left(
\begin{matrix} \left[\frac{1}{2}y_t^{(1)}-y_t^{(2)}\right]\Big{[}\left(y_t^{(1)}\right)^2+\left(y_t^{(2)}\right)^2\Big{]} \\ \left[y_t^{(1)}+\frac{1}{2}y_t^{(2)}\right]\Big{[}\left(y_t^{(1)}\right)^2+\left(y_t^{(2)}\right)^2\Big{]}\end{matrix}
\right)\mathrm{d}t-
\left(
\begin{matrix} z_t^{(1)} \\ z_t^{(2)} \end{matrix} \right)\mathrm{d}W_t, \\
&y_T=\left[
\begin{matrix}y_T^{(1)} \\ y_T^{(2)}\end{matrix} \right]=\left[
\begin{matrix}\sin\left(W_T+T\right) \\ \cos\left(W_T+T\right)\end{matrix} \right],\label{22:BSDE:gov}
\end{split}
\right.
\end{equation}
with the exact solutions£º
\begin{equation}
y_t=\left[
\begin{matrix}y_t^{(1)} \\ y_t^{(2)}\end{matrix} \right]=\left[
\begin{matrix}\sin\left(W_t+t\right) \\ \cos\left(W_t+t\right)\end{matrix} \right],\;\;\;\;
 z_t=\left[
\begin{matrix}z_t^{(1)} \\ z_t^{(2)}\end{matrix} \right]=\left[
\begin{matrix}\cos\left(W_t+t\right) \\ -\sin\left(W_t+t\right)\end{matrix} \right]. \label{22:analytic:solution}
\end{equation}
At the initial moment, the exact solutions  are: $y_0=[y_0^{(1)},y_0^{(2)}]^*=[0,1]^*$, $z_0=[z_0^{(1)},z_0^{(2)}]^*=[1,0]^*$. Without loss of generality, let us choose $T=1$. According to (\ref{2:bsde:corr:diff})-(\ref{2:bsde:solution:relation}), the corresponding PDEs read
\begin{equation}
\left\{
\begin{split}
{\cal N}_1[u_1,u_2]=\frac{\partial u_1}{\partial t}+\frac{1}{2}\frac{\partial ^2 u_1}{\partial x^2}+\left(\frac{1}{2}u_1-u_2\right)\Big{(} u_1^2+u_2^2\Big{)}&=0,\\
 {\cal N}_2[u_1,u_2]=\frac{\partial u_2}{\partial t}+\frac{1}{2}\frac{\partial ^2 u_2}{\partial x^2}+\left(u_1+\frac{1}{2}u_2 \right)\Big{(} u_1^2+u_2^2\Big{)}&=0,\label{22:pde:gov}
\end{split}
\right.
\end{equation}
subject to the boundary conditions
\begin{equation}
u_1(1,x)=\sin(x+1),\;\;\;\;\;\;\;\;u_2(1,x)=\cos(x+1), \label{22:pde:boundary}
\end{equation}
with the exact solutions
\begin{equation}
u_1(t,x)=\sin(x+t),\;\;\;\;\;\;\;\;u_2(t,x)=\cos(x+t).
\end{equation}
Here, ${\cal N}_1$ and ${\cal N}_2$ are nonlinear operators defined by (\ref{22:pde:gov}).  

In the frame of the HAM, we first of all construct a family of differential equations, called the zeroth-order deformation equations
\begin{eqnarray}\label{22:zeorth}
\left\{
\begin{split}
(1-q)\frac{\partial }{\partial t}\Big{[}\Phi(t,x,q)-\varphi_{0}(t,x)\Big{]}=c_0 \; q \; {\cal N}_1\Big{[}\Phi(t,x,q), \Xi(t,x,q) \Big{]}, \\
(1-q)\frac{\partial}{\partial t}\Big{[}\Xi(t,x,q)-S_{0}(t,x)\Big{]}=c_0 \; q \;  {\cal N}_2\Big{[}\Phi(t,x,q), \Xi(t,x,q) \Big{]}, 
\end{split}
\right.
\end{eqnarray}
subject to the boundary condition
\begin{equation}
\Phi(t,x,q)\Big{|}_{t=1}=\sin(x+1),\;\;\;\;\;\;\;\;\Xi(t,x,q)\Big{|}_{t=1}=\cos(x+1), \label{22:boundary}
\end{equation}
where $\varphi_0$ and $S_0$ are the initial guesses of $u_1$ and $u_2$, respectively.  It is obvious that
\begin{equation}
\Phi(t,x,1)=u_1(t,x),\;\;\;\;\;\;\;\;\Xi(t,x,1)=u_2(t,x).
\end{equation}
So, as $q$ increases from $0$ to $1$, $\Phi(t,x,q)$ varies continuously from its initial guess  $\varphi_0(t,x)$ to the unknown solution $u_1(t,x)$, $\Xi(t,x,q)$ varies continuously from its initial guess $S_0(t,x)$ to $u_2(t,x)$, respectively.

Similarly, we expand $\Phi(t,x,q)$ and $\Xi(t,x,q)$ into power series of $q$, i.e. 
\begin{equation}
\Phi(t,x,q)=\sum_{m=0}^{+\infty}\varphi_m(t,x) \cdot q^m,  \;\;\;\;\Xi(t,x,q)=\sum_{m=0}^{+\infty}S_m(t,x) \cdot q^m. \label{case2:series}
\end{equation}
Substituting (\ref{case2:series}) into Eqs.~(\ref{22:zeorth}) and (\ref{22:boundary}), then equating the like-power of $q$, we have the high-order deformation equations
\begin{eqnarray}
\frac{\partial}{\partial t}\Big{(}\varphi_m-\chi_m\varphi_{m-1}\Big{)}=c_0 \; \delta_{1,m-1},  \label{geq:mth:exam2:1st}\\
\frac{\partial}{\partial t}\Big{(}S_m-\chi_mS_{m-1}\Big{)}=c_0 \; \delta_{2,m-1},   \label{geq:mth:exam2:2nd}
\end{eqnarray}
subject to the boundary condition
\begin{equation}
\varphi_m(t,x)\Big{|}_{t=1}=0,\;\;\;\;\;\;\;\;S_m(t,x)\Big{|}_{t=1}=0,
\end{equation}
where $\chi_m$ is defined by (\ref{def:chi}), and 
\begin{eqnarray}
\delta_{1,n}(t,x) &=& {\cal D}_n \bigg{(} {\cal N}_1\Big{[}\Phi(t,x,q), \Xi(t,x,q) \Big{]} \bigg{)}\nonumber\\ 
&=& \frac{\partial \varphi_{n}}{\partial t}+\frac{1}{2}\frac{\partial^2 \varphi_{n}}{\partial x^2} \nonumber \\
&+&\sum_{i=0}^{n}\sum_{j=0}^{n-i}\left(\frac{1}{2}\varphi_i-S_i\right)\Big{(}\varphi_j\varphi_{n-i-j}+S_{j}S_{n-i-j}\Big{)}, \quad \quad \\
\delta_{2,n}(t,x) &=& {\cal D}_n \bigg{(} {\cal N}_2\Big{[}\Phi(t,x,q), \Xi(t,x,q) \Big{]} \bigg{)}\nonumber\\ 
&=& \frac{\partial S_{n}}{\partial t}+\frac{1}{2}\frac{\partial^2 S_{n}}{\partial x^2} \nonumber \\
&+&\sum_{i=0}^{n}\sum_{j=0}^{n-i}\left(\varphi_i+\frac{1}{2}S_i\right)\Big{(}\varphi_j\varphi_{n-i-j}+S_{j}S_{n-i-j}\Big{)},   
\end{eqnarray}
where the operator ${\cal D}_n$ is defined by (\ref{def:D}).   To satisfy the boundary condition (\ref{22:pde:boundary}), we choose the initial guesses
\begin{equation}
\varphi_0(t,x)=\sin(x+1),\;\;\;\;\;\;\;\;S_0(t,x)=\cos(x+1).
\end{equation}
Note that the high-order deformation equations (\ref{geq:mth:exam2:1st}) and (\ref{geq:mth:exam2:2nd}) are linear and {\em decoupled}, and thus can be easily obtained via integration with respect to $t$, say, 
\begin{eqnarray}
\varphi_m(t,x) &=& \chi_m \; \varphi_{m-1}(t,x) + c_0 \int_1^t \delta_{1,m-1}(t,x)\mathrm{d}t, \\
S_m(t,x) &=& \chi_m \; S_{m-1}(t,x) + c_0 \int_1^t \delta_{2,m-1}(t,x) \mathrm{d} t, 
\end{eqnarray}
step by step, starting from $m=1$, especially by means of computer algebra software such as Mathematica.   Similarly, the $M$th-order homotopy-approximations of $u_1(t,x)$ and $u_2(t,x)$ are given by
\begin{equation}
\tilde{\varphi}_M(t,x)=\sum_{m=0}^{M}\varphi_{m}(t,x),\;\;\;\;\;\;\;\;\tilde{S}_M(t,x)=\sum_{m=0}^{M}S_{m}(t,x),
\end{equation}
which give the approximations of the  initial values 
\begin{equation}
\begin{split}
\tilde{y}_1&=\tilde{\varphi}_M(0,x)\Big{|}_{x=0},\;\;\;\;\;\;\;\;\;\;\;\tilde{y}_2\;=\;\tilde{S}_M(0,x)\Big{|}_{x=0},\\
\tilde{z}_1&=\frac{\partial \tilde{\varphi}_M(0,x)}{\partial x}\bigg{|}_{x=0}, \;\;\;\;\;\;\;\;\tilde{z}_2\;=\;\frac{\partial \tilde{S}_M(0,x)}{\partial x}\bigg{|}_{x=0}, \nonumber
\end{split}
\end{equation}
respectively.   In this way, the original {\em coupled} nonlinear partial differential equations (\ref{22:pde:gov}) and (\ref{22:pde:boundary}) are solved just by means of integrations with respect to the time $t$.  

Define the squared residual errors
\begin{eqnarray}
\left\{
\begin{split}
\tilde{{\cal E}}&=\int_{0}^{T}\int_{-\infty}^{+\infty}\left[\Big{(}\tilde{\varphi}_M-u_1\Big{)}^2+\left(\tilde{S}_M-u_2\right)^2\right]\mathrm{d}x\mathrm{d}t,\\
{\cal E}&=\int_{0}^{T}\int_{-\infty}^{+\infty}\bigg{\{}\Big{(}{\cal N}_1[\tilde{\varphi}_M,\tilde{S}_M]\Big{)}^2+\Big{(}{\cal N}_2[\tilde{\varphi}_M,\tilde{S}_M]\Big{)}^2\bigg{\}}\mathrm{d}x\mathrm{d}t.
\end{split}
\right.
\end{eqnarray}
Similarly,  it is found that convergent results can be obtained with \emph{any} $c_0\in[-1.3,-0.4]$, but the optimal $c_0$ (i.e. the minimum $\cal E$) is near to $-1$. So we choose $c_0=-1$, and the homotopy approximations are
\begin{equation*}
\left\{
\begin{split}
\tilde{\varphi}_1(t,x)&=(t-1)\cos(1+x)+\sin(1+x);\\
\tilde{S}_1(t,x)&=(1-t)\sin(1+x)+\cos(1+x);\\
\tilde{\varphi}_2(t,x)&=(t-1)\cos(1+x)-\frac{1}{2}(t^2-2t-1)\sin(1+x);\\
\tilde{S}_2(t,x)&=(1-t)\sin(1+x)-\frac{1}{2}(t^2-2t-1)\cos(1+x);\\
\tilde{\varphi}_3(t,x)&=-\frac{1}{6}(t^3-3t^2-3t+5)\cos(1+x)-\frac{1}{2}(t^2-2t-1)\sin(1+x);\\
\tilde{S}_3(t,x)&=\frac{1}{6}(t^3-3t^2-3t+5)\sin(1+x)-\frac{1}{2}(t^2-2t-1)\cos(1+x);\\
\cdots\quad&
\end{split}
\right.
\end{equation*}

\begin{figure}[!t]
    \begin{center}
        \begin{tabular}{cc}
            \includegraphics[width=2.5in]{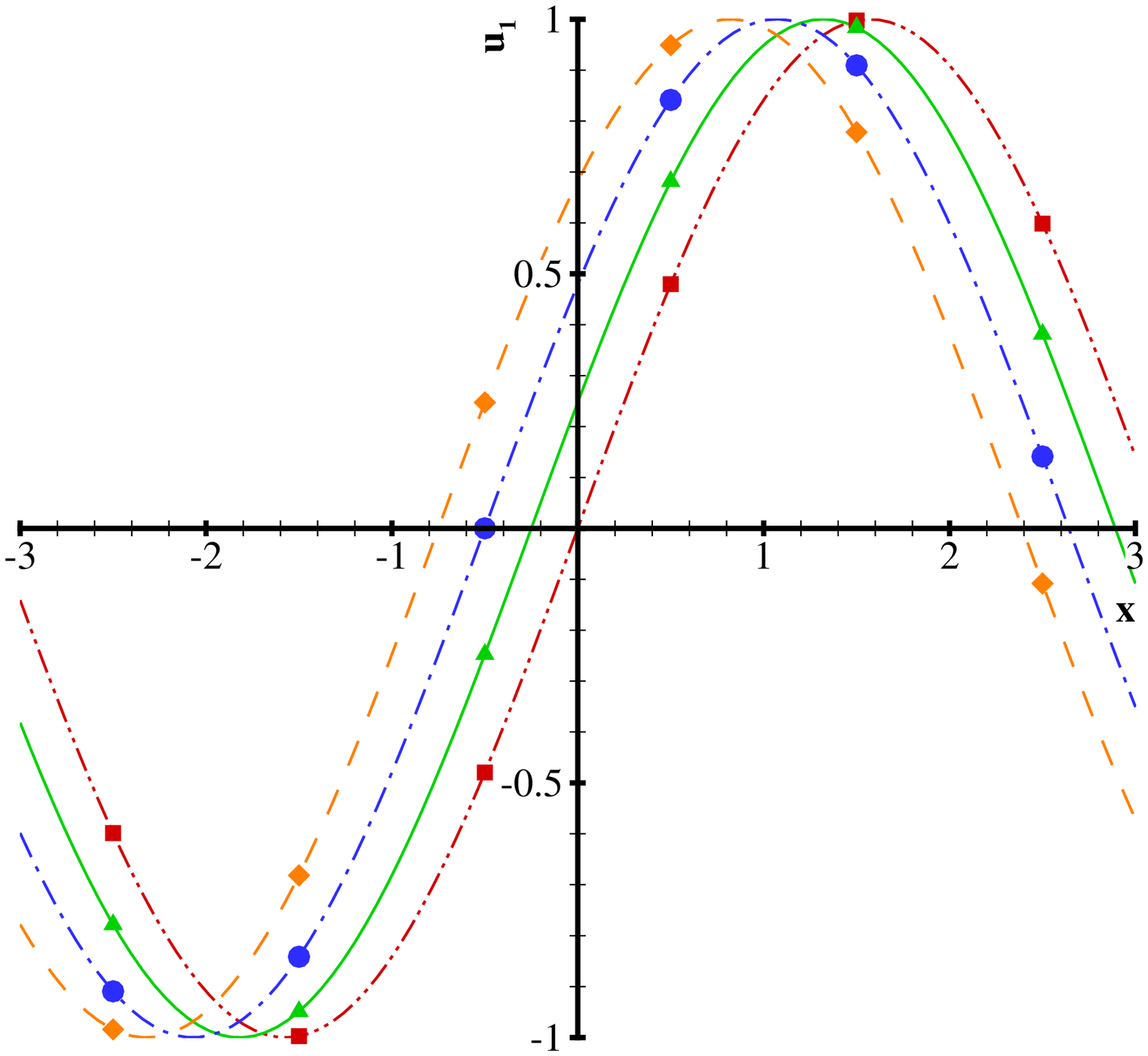}
            \includegraphics[width=2.5in]{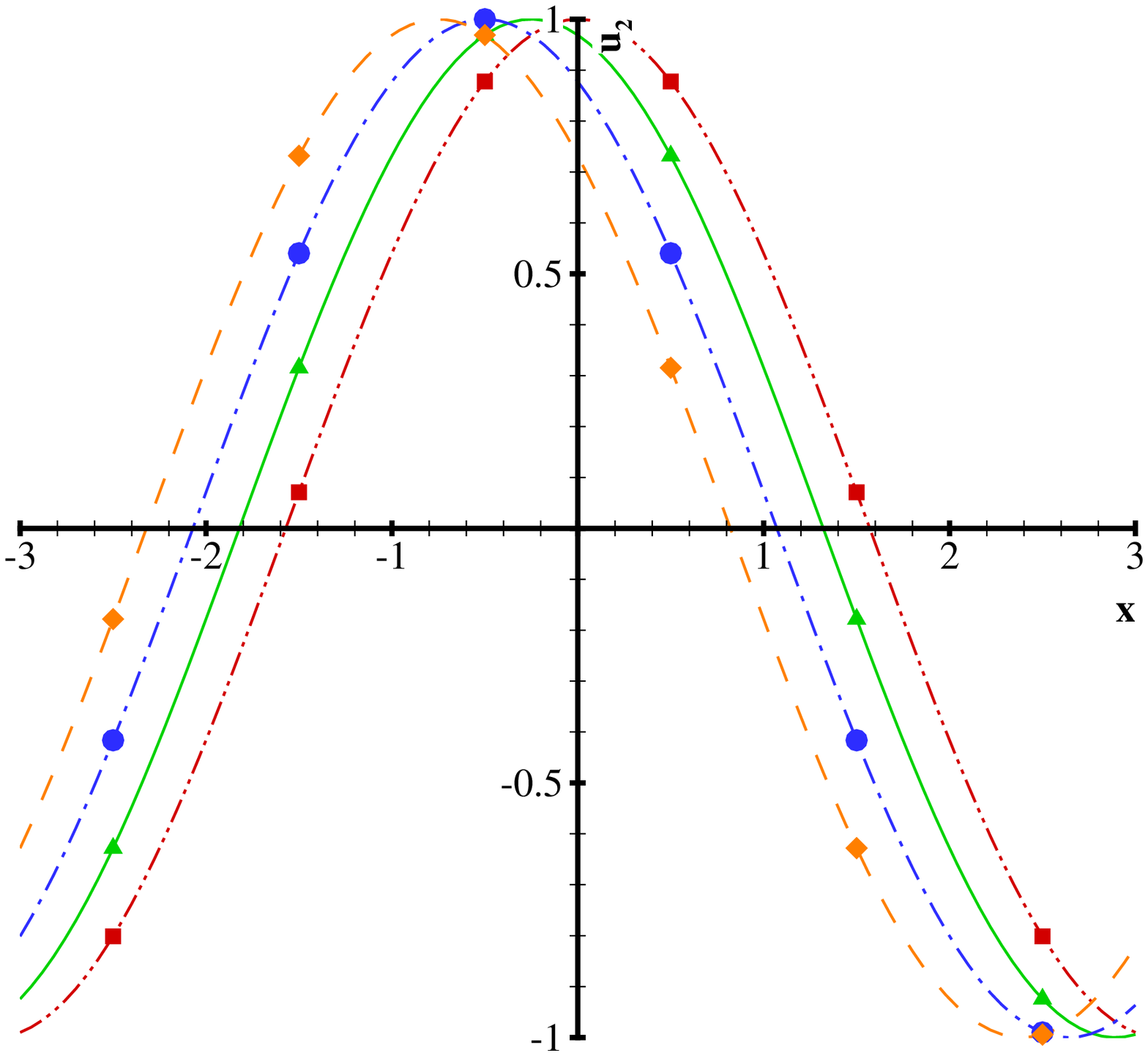}
        \end{tabular}
    \caption{Comparison between the homotopy approximations and the exact solutions $u_1$ (left), $u_2$ (right) of Eqs.~(\ref{22:pde:gov}) and (\ref{22:pde:boundary}), respectively, at $t=0, 0.25, 0.5, 0.75$. Dash-double-dotted line: homotopy approximations at $t=0$; Solid line: homotopy approximations at $t=0.25$; Dashed-dotted line: homotopy approximations at $t=0.5$; Dashed line: homotopy approximations at $t=0.75$. Square: the exact solution at $t=0$; Triangle up: the exact solution at $t=0.25$; Circle: the exact solution at $t=0.5$; Diamond: the exact solution at $t=0.75$.} \label{case2:u}
    \end{center}
\end{figure}

As shown in Table \ref{case2}, homotopy approximations with high accuracy at the level of $10^{-20}$ can be obtained in just 6.9 seconds CPU times.  In addition, our homotopy approximations agree well with the exact solutions, as shown in Fig.~\ref{case2:u}. So, this kind of BSDEs can also be solved by the HAM with high efficiency.
\begin{table}[!htb]
\tabcolsep 0pt
\caption{The squared residual error, the used CPU time (second) and the relative error between the homotopy approximations $\tilde{y}_1$, $\tilde{y}_2$ and the exact solutions $y_0^{(1)}$, $y_0^{(2)}$ of Eq.~(\ref{22:BSDE:gov}), respectively. Here, $m$ denotes the order of approximation.}
\vspace*{-12pt}\label{case2}
\begin{center}
\def\temptablewidth{1\textwidth}
{\rule{\temptablewidth}{1pt}}
\begin{tabular*}{\temptablewidth}{@{\extracolsep{\fill}}ccccc}
$m$    & $\tilde{\cal E}$& $\tilde{y}_1-y_{0}^{(1)}$ &$\tilde{y}_2-y_{0}^{(2)}$  & time (s)  \\
\hline
4 &   $8\times10^{-6}$ & $6\times10^{-3}$ & $-6\times10^{-3}$     &  0.2   \\
8 &    $6\times10^{-13}$ &$2\times10^{-6}$ & $-2\times10^{-6}$    &  0.7    \\
12 &   $2\times10^{-21}$ & $1\times10^{-10}$& $-1\times10^{-10}$   & 1.7 \\
16 &   $5\times10^{-31}$ & $2\times10^{-15}$& $-2\times10^{-15}$   & 3.5 \\
20 &   $2\times10^{-41}$ & $1\times10^{-20}$& $-2\times10^{-20}$     &  6.9\\
\end{tabular*}
{\rule{\temptablewidth}{1pt}}
 \end{center}
 \end{table}

 \subsection{An example of BSDE with two-dimensional Brownian motion} %
 
 In this case, we consider the following BSDEs with two dimensional Brownian motion
 \begin{eqnarray}
\left\{
\begin{split}
-\mathrm{d}y_t=\left(y_t-\frac{1}{2}z_t^{(1)}-\frac{1}{2}z_t^{(2)}\right)\mathrm{d}t-z_t^{(1)}\mathrm{d}W_t^{(1)}-z_t^{(2)}\mathrm{d}W_t^{(2)}, \;0\leq t\leq T,\;\;\\
y_T=\sin\left(W_T^{(1)}+W_T^{(2)}+T\right),\quad\quad\quad\quad\quad\quad\quad~~~~~~~~~~~~~~~~~~~~~~~~~~~~~~~ \label{23:case2:bsde}
\end{split}
\right.
\end{eqnarray}
with the exact solutions
\begin{eqnarray}
\left\{
\begin{split}
y_t&=&\sin\left(W_t^{(1)}+W_t^{(2)}+t\right),\\
z_t^{(1)}&=&\cos\left(W_t^{(1)}+W_t^{(2)}+t\right),\\
z_t^{(2)}&=&\cos\left(W_t^{(1)}+W_t^{(2)}+t\right).\label{21:case1:solution}
\end{split}
\right.
\end{eqnarray}
At the initial moment, the exact solutions are: $y_0=0$, $z_0^{(1)}=1$, $z_0^{(2)}=1$. Without loss of generality, we choose $T=1$. According to (\ref{2:bsde:corr:diff})-(\ref{2:bsde:solution:relation}), the corresponding PDE reads
\begin{equation}
{\cal N}_3[u]=\frac{\partial u}{\partial t}+\frac{1}{2}\frac{\partial ^2 u}{\partial x_1^2}+\frac{1}{2}\frac{\partial ^2 u}{\partial x_2^2}+u-\frac{1}{2}\frac{\partial u}{\partial x_1}-\frac{1}{2}\frac{\partial u}{\partial x_2}=0,  \label{23:pde:gov}
\end{equation}
subjects to the boundary condition
\begin{equation}
u(1,x_1,x_2)=\sin(x_1+x_2+1). \label{23:pde:boundary}
\end{equation}
Here ${\cal N}_3$ is a nonlinear operator, defined by (\ref{23:pde:gov}).  

Similarly, in the frame of the HAM, we construct the following zeroth-order deformation equations
\begin{equation}
(1-q)\; \frac{\partial }{\partial t}\Big{[}\Phi(t,x_1,x_2,q)-\varphi_{0}(t,x_1,x_2)\Big{]} =c_0 \; q \; {\cal N}_3 \Big{[}\Phi(t,x_1,x_2,q) \Big{]}, \label{23:zeroth}
\end{equation}
subject to the boundary condition
\begin{equation}
\Phi(t,x_1,x_2,q)\Big{|}_{t=1}=\sin(x_1+x_2+1),
\end{equation}
where $\varphi_0$ is an initial guess of the unknown $u(t,x_1,x_2)$, which satisfies the boundary condition (\ref{23:pde:boundary}).    The above equations construct a continuous variation $\Phi(t,x_1,x_2,q)$ from the initial guess $\varphi_0(t,x_1,x_2)$ at $q=0$ to the unknown $u(t,x_1,x_2)$ at $q=1$.

Similarly, the $M$th-order homotopy-approximation of $u(t,x_1,x_2)$ is given by
\begin{equation}
\tilde{\varphi}_M(t,x_1,x_2)=\sum_{m=0}^{M}\varphi_{m}(t,x_1,x_2),
\end{equation}
where $\varphi_{m}$ is governed by the high-order deformation equations
\begin{equation}
\begin{split}
\frac{\partial}{\partial t}\Big{(}\varphi_m-\chi_m\varphi_{m-1}\Big{)}&=c_0 \; \delta_{3,m-1},
\end{split}
\end{equation}
subject to the boundary condition
\begin{equation}
\varphi_m(t,x_1,x_2)\Big{|}_{t=1}=0,
\end{equation}
where 
\begin{eqnarray}
&&\delta_{3,k}(t,x_1,x_2) = {\cal D}_k \bigg{(}  {\cal N}_3 \Big{[}\Phi(t,x_1,x_2,q) \Big{]} \bigg{)}\nonumber\\
&=& \frac{\partial \varphi_{k}}{\partial t}+\frac{1}{2}\left( \frac{\partial^2 \varphi_{k}}{\partial x_1^2}+\frac{\partial^2 \varphi_{k}}{\partial x_2^2}\right) +\varphi_{k} -\frac{1}{2}\left(\frac{\partial \varphi_{k}}{\partial x_1}+\frac{\partial \varphi_{k}}{\partial x_2}\right), 
\end{eqnarray}
and $\chi_m$ is defined by (\ref{def:chi}),  ${\cal D}_k$ is defined by (\ref{def:D}), respectively.    
According to (\ref{23:pde:boundary}), we choose the initial guess
\begin{equation}
\varphi_0(t,x_1,x_2)=\sin(x_1+x_2+1).
\end{equation}
Note that it is straightforward to gain the solution 
\begin{equation}
\varphi_m(t,x_1,x_2) = \chi_m \varphi_{m-1}(t,x_1,x_2) + c_0 \; \int_1^t \delta_{3,m-1}(z,x_1,x_2) \mathrm{d} z, 
\end{equation}
step by step, starting from $m=1$, especially by means of the computer algebra software Mathematica.   In other words,  in the frame of the HAM, the original partial differential equations (\ref{23:pde:gov}) and (\ref{23:pde:boundary}) are solved just by means of integrations with respect to the time $t$.   As long as the $M$th-order approximation $\tilde{\varphi}_M$ is obtained,  we have the approximations of the initial values  
\begin{equation}
\tilde{y}_0=\tilde{\varphi}_M(0,0,0),\;\;\;\tilde{z}_0^{(1)}=\frac{\partial \tilde{\varphi}_M(0,x_1,0)}{\partial x_1}\Big{|}_{x_1=0},\;\;\;\tilde{z}_0^{(2)}=\frac{\partial \tilde{\varphi}_M(0,0,x_2)}{\partial x_2}\Big{|}_{x_2=0}. \nonumber
\end{equation} 
Thanks to the symmetry, it is obvious that
\begin{equation}
\tilde{z}_0^{(1)}-z_0^{(1)}=\tilde{z}_0^{(2)}-z_0^{(2)}.
\end{equation}

\begin{figure}[!h]
    \begin{center}
        \begin{tabular}{cc}
            \includegraphics[width=3in]{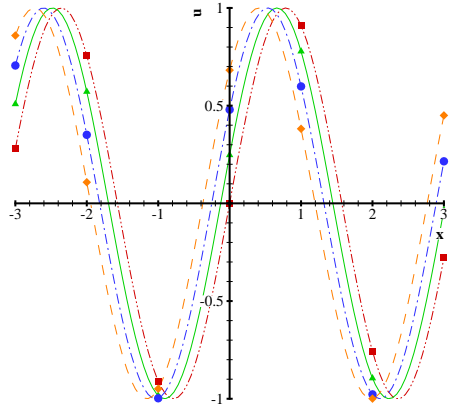}
        \end{tabular}
    \caption{Comparison between the homotopy approximation and the exact solution $u$ of Eqs.~(\ref{23:pde:gov}) and (\ref{23:pde:boundary}), when $x_1=x_2=x$, $t=0, 0.25, 0.5, 0.75$. Dash-double-dotted line: homotopy approximation at $t=0$; Solid line: homotopy approximation at $t=0.25$; Dashed-dotted line: homotopy approximation at $t=0.5$; Dashed line: homotopy approximation at $t=0.75$. Square: exact solution at $t=0$; Triangle up: exact solution at $t=0.25$; Circle: exact solution at $t=0.5$; Diamond: exact solution at $t=0.75$.} \label{case3:u}
    \end{center}
\end{figure}

Similarly, define the squared residual errors
\begin{equation}
\left\{
\begin{split}
\tilde{{\cal E}}&=\int_{0}^{T}\int_{-\infty}^{+\infty}\int_{-\infty}^{+\infty}\Big{[}\left(\tilde{\varphi}_M-u\right)^2\Big{]}\mathrm{d}x_1\mathrm{d}x_2\mathrm{d}t,\\
{\cal E}&=\int_{0}^{T}\int_{-\infty}^{+\infty}\int_{-\infty}^{+\infty}\Big{\{}{\cal N}_3[\tilde{\varphi}_M(t,x_1,x_2)]\Big{\}}^2\mathrm{d}x_1\mathrm{d}x_2\mathrm{d}t.
\end{split}
\right.
\end{equation}
It is found that convergent results can be obtained with \emph{any} $c_0\in[-1.3,-0.4]$, but the optimal $c_0$ (i.e. the minimum $\cal E$) is near to $-1$. So,  we choose $c_0=-1$, and the corresponding homotopy approximations are
\begin{equation*}
\left\{
\begin{split}
\tilde{\varphi}_1(t,x)&=(t-1)\cos(1+x_1+x_2)+\sin(1+x_1+x_2);\\
\tilde{\varphi}_2(t,x)&=(t-1)\cos(1+x_1+x_2)-\frac{1}{2}(t^2-2t-1)\sin(1+x_1+x_2);\\
\tilde{\varphi}_3(t,x)&=-\frac{1}{6}(t^3-3t^2-3t+5)\cos(1+x_1+x_2)\\
&\quad-\frac{1}{2}(t^2-2t-1)\sin(1+x_1+x_2);\\
\cdots\quad&
\end{split}
\right.
\end{equation*}
As shown in Table \ref{case3}, accurate approximations at the level of $10^{-17}$ are obtained in just 0.5 seconds CPU times.  Besides, our homotopy approximations agree well with the exact solution, as shown in Fig.~\ref{case3:u}. This again demonstrates the validity and high efficiency of the HAM for the BSDEs.
\begin{table}[!htb]
\tabcolsep 0pt
\caption{The squared residual error, used CPU time (second) and the relative error between the homotopy approximations $\tilde{y}_0$, $\tilde{z}_0^{(1)}$ and the exact solutions $ y_0$, $z_0^{(1)}$ of Eq.~(\ref{23:case2:bsde}), respectively. Here, $m$ denotes the order of approximation.}
\vspace*{-12pt}\label{case3}
\begin{center}
\def\temptablewidth{1\textwidth}
{\rule{\temptablewidth}{1pt}}
\begin{tabular*}{\temptablewidth}{@{\extracolsep{\fill}}ccccc}
$m$    &  $\tilde{\cal E}$ & $\tilde{y}_0-y_{0}$ &$\tilde{z}_0^{(1)}-z_{0}^{(1)}$     & time (s)  \\
\hline
4 &  $4\times10^{-6}$ &  $6\times10^{-3}$ & $-6\times10^{-3}$    &   0.06   \\
8 &   $3\times10^{-13}$ & $2\times10^{-6}$ & $-2\times10^{-6}$  &   0.12    \\
12 &  $8\times10^{-22}$ &  $1\times10^{-10}$& $-1\times10^{-10}$  &  0.23 \\
16 &  $2\times10^{-31}$ &  $2\times10^{-15}$& $-2\times10^{-15}$   & 0.36 \\
20 & $1\times10^{-41}$ &   $-4\times10^{-17}$& $7\times10^{-17}$ &   0.50\\
\end{tabular*}
{\rule{\temptablewidth}{1pt}}
 \end{center}
 \end{table}

Here, it should be emphasized that, although the three kinds of BSDEs are quite different,  our HAM approaches have no fundamental difference: the convergent approximations are always gained  simply by means of integration with respect to the time $t$, within a few seconds CPU times!         

\renewcommand{\theequation}{3.\arabic{equation}}
\setcounter{equation}{0}
\section{The HAM for FBSDEs}%

In this section, we further employ the HAM to solve some FBSDEs, 2FBSDEs and high-dimensional FBSDEs.

\subsection{An example of FBSDEs} %

First, let us consider the following coupled FBSDE
\begin{eqnarray}
\left\{
\begin{split}
\mathrm{d}x_t&=&y_t\sin(t+x_t)\mathrm{d}t+y_t\cos(t+x_t)\mathrm{d}W_t,\quad\quad\quad\quad\quad\quad\quad~~~~\\
-\mathrm{d}y_t&=&\left[\frac{1}{2}y_t^2z_t-\cos(t+x_t)(y_t^2+1) \right]\mathrm{d}t-z_t\mathrm{d}W_t, \;\;\;\;0\leq t\leq T,\\
x_0&=&\frac{1}{2},\;\;\;\;\;\;\;\;y_T=\sin(T+x_T),\quad\quad\quad~~~~~~~~~~~~~~~~~~~~~~~~~~~ \label{3:fbsde}
\end{split}
\right.
\end{eqnarray}
with the exact solutions
\begin{equation}
y_t=\sin(t+x_t),\quad\quad\quad
z_t=\cos^2(t+x_t)\sin(t+x_t).\label{3:analytic:solution}
\end{equation}
At the initial moment, the exact values read
 \begin{equation}
 y_0=\sin\left(x_0\right),\quad\quad z_0=\cos^2\left(x_0\right)\sin\left(x_0\right). \nonumber
 \end{equation}
Without loss of generality, we choose $T=1$. Then according to Eqs.~(\ref{1:general:fbsde})-(\ref{1:general:fbsde:boundary}), the corresponding PDE reads
\begin{eqnarray}
\frac{\partial u}{\partial t} &+&\sin(t+x)u\frac{\partial u}{\partial x}+\frac{1}{4}\cos(2t+2x)u^2\frac{\partial^2 u}{\partial x^2}+\frac{1}{4}u^2\frac{\partial^2 u}{\partial x^2} \nonumber \\
&+&\frac{1}{2}\cos(t+x)u^3 \frac{\partial u}{\partial x}-\cos(t+x)(u^2+1)=0, \label{31:origin:equation}
\end{eqnarray}
subjects to the boundary condition
\begin{equation}
u(t,x)\Big{|}_{t=1}=\sin(1+x), \label{31:pde:boundary}
\end{equation}
with the exact solution $u=\sin(t+x)$.     

In the frame of the HAM, we first construct the zeroth-order deformation equations
\begin{eqnarray}
\begin{split}
&(1-q)\frac{\partial }{\partial t}\Big{[}\Phi(t,x,q)-\varphi_{0}(t,x)\Big{]}=c_0\; q \; {\cal N}_4\Big{[}\Phi(t,x,q),q \Big{]}, \label{31:zeroth}
\end{split}
\end{eqnarray}
subject to the boundary condition
\begin{equation}
\Phi(t,x,q)\Big{|}_{t=1}=\sin(q+x), \label{31:zeroth:boundary}
\end{equation}
where $\varphi_0$ is an initial guess of $u(t,x)$, and 
\begin{eqnarray}
&&{\cal N}_4\Big{[}\Phi(t,x,q),q \Big{]} \nonumber\\
&=&\frac{\partial \Phi}{\partial t} +\sin(t q+x)u\frac{\partial \Phi}{\partial x}+\frac{1}{4}\cos(2tq+2x)\Phi^2\frac{\partial^2 \Phi}{\partial x^2}+\frac{1}{4}\Phi^2\frac{\partial^2 \Phi}{\partial x^2} \nonumber \\
&+&\frac{1}{2}\cos(tq+x)\Phi^3 \frac{\partial \Phi}{\partial x}-\cos(tq+x)(\Phi^2+1) \label{def:N:4}
\end{eqnarray}
corresponds to the original nonlinear differential equation (\ref{31:origin:equation}), 
respectively.   Note that, in order to express $u(t,x)$ as a power series of $t$,  we replace $\sin(t+x)$, $\cos(t+x)$  and $\cos(2t+2x)$  by $\sin(t q+x)$, $\cos(t q+x)$ and $\cos(2t q+2x)$ in the governing equation  (\ref{31:zeroth}), respectively, where $q\in[0,1]$ is the embedding parameter.   Besides, we replace $\sin(1+x)$  in the boundary condition (\ref{31:zeroth:boundary}) by $\sin(q+x)$.   
Otherwise,   the solution $\varphi_m$ contains the complicated terms such as $t^n\cos(t+x)$ and $t^n \sin(t+x)$, which do not belong to the power series of $t$  and thus  lead to difficulty to gain high-order approximations via integration with respect to the time $t$.      It should be emphasized that,   
 we can do all of these  mainly because the HAM provides us such kind of great freedom to construct zeroth-order deformation equations \cite{Liaobook, liaobook2}.   
 
 Similarly,  the above equations (\ref{31:zeroth}) and (\ref{31:zeroth:boundary}) construct a continuous deformation (i.e. a homotopy) $\Phi(t,x,q)$ from the initial guess $\varphi_0$ to the unknown solution $u(t,x)$ as $q$ increases from 0 to 1.   Similarly, we can expand $\Phi(t,x,q)$ into a power series of $q$, i.e.
\begin{equation}
\Phi(t,x,q)=\sum_{m=0}^{+\infty}\varphi_{m}(t,x)\cdot q^m. \label{31:power:series}
\end{equation}
Differentiating the zeroth-order deformation equations (\ref{31:zeroth}) and (\ref{31:zeroth:boundary}) $m$ times with respect to $q$, then setting $q=0$, and finally dividing the two sides of them by $m!$, we have the high-order deformation equations
\begin{eqnarray}
\frac{\partial}{\partial t}\Big{(}\varphi_m-\chi_m\varphi_{m-1}\Big{)} = c_0\; \delta_{4,m-1}, \label{31:highorder:delta}
\end{eqnarray}
subject to the boundary condition
\begin{equation}
\varphi_m(t,x)\Big{|}_{t=1}=\frac{1}{m!}\sin\left(x+\frac{m \pi}{2}\right). \label{31:highorder:boundary}
\end{equation}
where 
\begin{eqnarray}
\delta_{4,n}(t,x) &=& {\cal D}_n \bigg{(} {\cal N}_4\Big{[}\Phi(t,x,q),q \Big{]} \bigg{)} \nonumber\\
& = & 
\frac{\partial \varphi_{n}}{\partial t}+\sum_{i=0}^{n}\sum_{j=0}^{i}\frac{t^{j}}{j!}\sin\left(x+\frac{j \pi }{2}\right)\varphi_{n-i}\frac{\partial \varphi_{i-j}}{\partial x} \nonumber \\
&+&\sum_{i=0}^{n}\sum_{j=0}^{i}\sum_{k=0}^{j}\sum_{p=0}^{k}\frac{t^p}{2\cdot p!}\cos\left(x+\frac{p \pi}{2}\right)\varphi_{n-i}\varphi_{i-j}\varphi_{j-k}\frac{\partial \varphi_{k-p}}{\partial x}\nonumber \\
&+&\sum_{i=0}^{n}\sum_{j=0}^{i}\sum_{k=0}^{j}\frac{(2t)^{k}}{4\cdot k!}\cos\left(2x+\frac{k \pi }{2}\right)\varphi_{n-i}\varphi_{i-j}\frac{\partial^2\varphi_{j-k}}{\partial x^2}\nonumber \\
&+&\sum_{i=0}^{n}\sum_{j=0}^{i}\left[\frac{1}{4}\frac{\partial^2\varphi_{j}}{\partial x^2}-\frac{t^j}{j!}\cos\left(x+\frac{j \pi}{2}\right)\right]\varphi_{n-i}\varphi_{i-j} \nonumber \\
&-&\frac{t^{n}}{n!}\cos\left(x+\frac{n\pi}{2}\right),   
\end{eqnarray}
and $\chi_m$ is defined by (\ref{def:chi}),  ${\cal D}_n$ is defined by (\ref{def:D}), ${\cal N}_4$ is defined by (\ref{def:N:4}),  respectively.    
According to (\ref{31:zeroth:boundary}), the initial guess of $u(t,x)$ should be
\begin{equation}
\varphi_0(t,x)=\sin(x).  \label{31:initial:guess}
\end{equation}
Then, it is easy to gain the solution 
\begin{equation}
\varphi_m(t,x) = \chi_m \; \varphi_{m-1}(t,x) + c_0 \int_0^t \delta_{4,m-1}(t,x)\mathrm{d} t + A_m(x), 
\end{equation}
where the integration coefficient $A_m(x)$ is determined by the boundary condition (\ref{31:highorder:boundary}), step by step, starting from $m=1$.   

The $M$th-order homotopy-approximation of $u(t,x)$ is gained at $q=1$ by
\begin{equation}
\tilde{\varphi}_M(t,x)=\sum_{m=0}^{M}\varphi_{m}(t,x), \label{31:homotopy:appro}
\end{equation}
which gives the approximations of the initial values
\begin{equation}
\tilde{y}_0=\tilde{\varphi}_M(0,x_0),\;\;\;\;\;\;\;\tilde{z}_0=\bigg{[}\cos(x)\tilde{\varphi}_M\frac{\partial \tilde{\varphi}_M}{\partial x}\bigg{]}\bigg{|}_{t=0,\;x=x_0}. \nonumber
\end{equation}

\begin{figure}[!t]
    \begin{center}
        \begin{tabular}{cc}
            \includegraphics[width=3in]{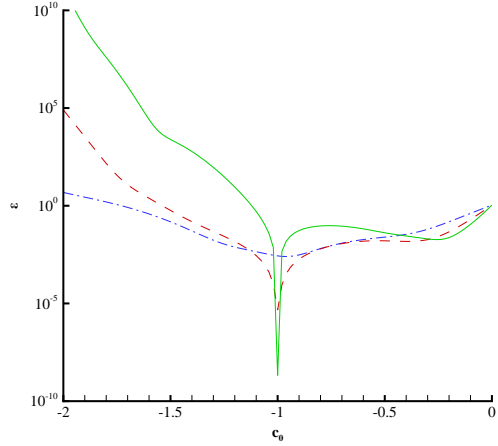}
        \end{tabular}
    \caption{The squared residual error versus $c_0$ of Eq.~(\ref{31:zeroth}). Dash-dotted line: 3th-order approximation; Dashed line: 5th-order approximation; Solid line: 7th-order approximation.} \label{case4:c0}
    \end{center}
\end{figure}

Similarly, the curves of $\cal E$ versus $c_0$ is as shown in Fig.~\ref{case4:c0}, which indicates that the optimal $c_0$ is $-1$. It is found that the homotopy-series is convergent within a rather small region. So, we choose the optimal $c_0=-1$, and the corresponding homotopy approximations are:
\begin{eqnarray}
\left\{
\begin{split}
\tilde{\varphi}_1(t,x)&=\sin(x)+t\cos(x);\\
\tilde{\varphi}_2(t,x)&=\sin(x)+t\cos(x)-\frac{1}{2}t^2\sin(x);\\
\tilde{\varphi}_3(t,x)&=\sin(x)+t\cos(x)-\frac{1}{2}t^2\sin(x)-\frac{1}{6}t^3\cos(x);\\
\cdots\quad&
\end{split}
\right.
\end{eqnarray}

It is interesting that the homotopy approximation $\tilde{\varphi}_M(t,x)$ is actually the Maclaurin series of the exact solution $u(t,x)=\sin(x+t)$ about $t$, say
\begin{equation}
\tilde{\varphi}_M(t,x)=\sum_{i=0}^{M}\frac{t^i}{i!}\sin\left(x+\frac{i}{2}\pi\right). \label{31:Maclaurin}
\end{equation}
A proof is given in the Appendix.  This also means that our homotopy-series converges to the exact solution of the original equations. Besides, our homotopy approximations at the initial moment $\tilde{y}_0$ and $\tilde{z}_0$ are exactly the same as the exact solutions $ y_0$ and $z_0$, respectively, and the $\tilde{\cal E}$ which characterizes the global error between the homotopy approximation and the exact solution quickly decreases to $9\times 10^{-42}$ within only 1.81 seconds CPU times, as shown in Table \ref{case4}.  It should be emphasized that the HAM provides us great freedom to  construct such kind of zeroth-order deformation equations (\ref{31:zeroth}) and (\ref{31:zeroth:boundary}).   In fact,  it is such kinds of freedom that distinguishes the HAM from other analytic methods.  All  of  these demonstrate the power and high efficiency of the HAM for forward-backward nonlinear differential equations.

\begin{table}[!htb]
\tabcolsep 0pt
\caption{The squared residual error, the used CPU time (second) and the relative error between the homotopy approximations $\tilde{y}_0$, $\tilde{z}_0$ and the exact solutions $ y_0$, $z_0$ of Eq.~(\ref{3:fbsde}), respectively. Here, $m$ denotes the order of approximation.}
\vspace*{-12pt}\label{case4}
\begin{center}
\def\temptablewidth{1\textwidth}
{\rule{\temptablewidth}{1pt}}
\begin{tabular*}{\temptablewidth}{@{\extracolsep{\fill}}ccccc}
$m$  &$\tilde{\cal E}$&  $\tilde{y}_0-y_{0}$ &$\tilde{z}_0-z_{0}$   & CPU time (s)  \\
\hline
4   & $4\times 10^{-6}$&  0 & 0    &   0.02   \\
8 & $3\times 10^{-13}$  & 0 & 0  &   0.11    \\
12  & $7\times 10^{-22}$&  0& 0 &  0.33 \\
16  & $2\times 10^{-31}$& 0 & 0   & 0.78 \\
20  & $9\times 10^{-42}$& 0 & 0 &   1.81\\
\end{tabular*}
{\rule{\temptablewidth}{1pt}}
 \end{center}
 \end{table}

\subsection{An example of 2FBSDEs} %
Then, let us consider the following 2FBSDE
\begin{equation}
\left\{
\begin{split}
\mathrm{d}x_t&=\sin(t+x_t)\mathrm{d}t+\cos(t+x_t)\mathrm{d}W_t,\\
-\mathrm{d}y_t&=\left[-\cos(t+x_t)z_t-\cos(t+x_t)(y_t^2+y_t)-\frac{\Gamma_t}{4} \right]\mathrm{d}t-z_t\mathrm{d}W_t,\\
\mathrm{d}z_t&=A_t\mathrm{d}t+\Gamma_t\mathrm{d}W_t,\;\;\;\;\quad\quad0\leq t\leq T,\\
x_0&=\frac{1}{2},\;\;\;\;\;\;\;\;y_T=\sin(T+x_T), \label{4:fbsde}
\end{split}
\right.
\end{equation}
with the exact solutions
\begin{equation}
\left\{
\begin{split}
y_t&=\sin(t+x_t),\quad \quad z_t=\cos^2(t+x_t), \\
\Gamma_t&=-2\sin(t+x_t)\cos^2(t+x_t),\\
A_t&=-\sin(2t+2x_t)\left[1+\sin(t+x_t)\right]-\cos(2t+2x_t)\cos^2(t+x_t).\label{4:analytic:solution}
\end{split}
\right.
\end{equation}
At the initial moment, the exact solutions  are
\begin{equation}
\left\{
\begin{split}
y_0&=\sin\left(x_0\right),\quad\quad z_0=\cos^2\left(x_0\right),\quad \quad \Gamma_0=-2\sin(x_0)\cos^2(x_0), \nonumber\\
A_0&=-\sin(2x_0)\left[1+\sin(x_0)\right]-\cos(2x_0)\cos^2(x_0).
\end{split}
\right.
\end{equation}
Without loss of generality, let us choose $T=1$. Then according to Eqs.~(\ref{1:general:2fbsde})-(\ref{ellip:ope}), the corresponding PDE reads
\begin{eqnarray}
&&\frac{\partial u}{\partial t}+\frac{1}{8}\left[1+\cos(2t+2x)\right]\frac{\partial^2 u}{\partial x^2}-\cos(t+x)(u^2+u) \nonumber \\
&+&\left[\sin(t+x)+\frac{1}{8}\sin(2t+2x)-\frac{1}{2}\cos(2t+2x)-\frac{1}{2} \right]\frac{\partial u}{\partial x} =0, \hspace{1.0cm} \label{4:origin:equation}
\end{eqnarray}
subjects to the boundary condition
\begin{equation}
u(t,x)\Big{|}_{t=1}=\sin(1+x). \label{32:pde:boundary}
\end{equation}

Similar to $\S~3.1$, we construct the zeroth-order deformation equations
\begin{equation}
(1-q)\frac{\partial }{\partial t}\Big{[}\Phi(t,x,q)-\varphi_{0}(t,x)\Big{]} = c_0\; q\; {\cal N}_5 \Big{[} \Phi(t,x,q),q \Big{]}
\end{equation}
subject to the boundary condition
\begin{equation}
\Phi(t,x,q)\Big{|}_{t=1}=\sin(q+x),
\end{equation}
where
\begin{eqnarray}
&&{\cal N}_5 \Big{[} \Phi(t,x,q),q \Big{]} \nonumber\\
&=&\frac{\partial \Phi}{\partial t}+\frac{1}{8}\left[1+\cos(2tq+2x)\right]\frac{\partial^2 \Phi}{\partial x^2}-\cos(tq+x)(\Phi^2+\Phi) \nonumber \\
&+&\left[\sin(tq+x)+\frac{1}{8}\sin(2tq+2x)-\frac{1}{2}\cos(2tq+2x)-\frac{1}{2} \right]\frac{\partial \Phi}{\partial x} \hspace{1.0cm}
\label{def:N:5}
\end{eqnarray}
corresponds to the original nonlinear equation (\ref{4:origin:equation}).  Note that we replace the terms $\cos(t+x), \sin(t+x)$ and so on by $\cos(tq+x)$ and $\sin(t q +x)$,  respectively, so as to gain a solution in the power series of $t$.

Similarly, the $M$th-order approximation is given at $q=1$ by (\ref{31:homotopy:appro}), with $\varphi_m$ being governed by the high-order deformation equations
\begin{eqnarray}
\frac{\partial}{\partial t}\Big{(}\varphi_m-\chi_m\varphi_{m-1}\Big{)} = c_0 \; \delta_{5,m-1},  \label{geq:exam5:mth}
\end{eqnarray}
subject to the boundary condition
\begin{equation}
\varphi_m(t,x)\Big{|}_{t=1}=\frac{1}{m!}\sin\left(x+\frac{m \pi}{2}\right),   \label{bc:exam5:mth}
\end{equation}
where
\begin{eqnarray}
\delta_{5,n} &=& {\cal D}_{n} \bigg{(} {\cal N}_5 \Big{[} \Phi(t,x,q),q \Big{]}  \bigg{)} \nonumber\\
&=&\frac{\partial \varphi_{n}}{\partial t}-\frac{1}{2}\frac{\partial \varphi_{n}}{\partial x}+\frac{1}{8}\frac{\partial^2 \varphi_{n}}{\partial x^2}+\sum_{i=0}^{n}\frac{t^{i}}{i!}\sin\left(x+\frac{i\pi}{2}\right)\frac{\partial \varphi_{n-i}}{\partial x} \nonumber\\
&+&\sum_{i=0}^{n}\frac{t^i}{i!}\left[\frac{2^i}{8}\cos\left(2x+\frac{i\pi}{2}\right)\frac{\partial^2 \varphi_{n-i}}{\partial x^2}-  \cos\left(x+\frac{i\pi}{2}\right)\varphi_{n-i} \right]  \nonumber\\
&+&\sum_{i=0}^{n}\frac{(2t)^i}{i!}\bigg{[}\frac{1}{8}\sin\left(2x+\frac{i\pi}{2}\right)-\frac{1}{2} \cos\left(2x+\frac{i\pi}{2}\right) \bigg{]}\frac{\partial \varphi_{n-i}}{\partial x} \nonumber\\
&-&\sum_{i=0}^{n}\sum_{j=0}^{n-i}\frac{t^i}{i!}\cos\left(x+\frac{i\pi}{2}\right)\varphi_j\varphi_{n-i-j},
\end{eqnarray}
with $\chi_m$, ${\cal D}_n$ and ${\cal N}_5$ being defined by (\ref{def:chi}), (\ref{def:D}) and (\ref{def:N:5}), respectively.    Similarly,  we choose $\varphi_0=\sin(x)$ as the initial guess of $u(t,x)$.   Then, it is easy to gain the solution of (\ref{geq:exam5:mth}): 
\begin{equation}
\varphi_m(t,x) = \chi_m \; \varphi_{m-1}(t,x) + c_0 \; \int_0^t \delta_{5,m-1}(t,x) \mathrm{d} t + B_m(x), 
\end{equation}
step by step, starting from $m=1$, where the integration coefficient $B_m(x)$ is determined by (\ref{bc:exam5:mth}).   This is  rather  efficient computationally by means of the computer algebra software such as Mathematica.  

As long as the $M$th-order approximation $\tilde{\varphi}_M$ is known,  we have the approximations of the initial values
\begin{eqnarray}
\left\{
 \begin{split}
  \tilde{y}_0&=\tilde{\varphi}_M(0,x_0),\;\;\;\;\;\;\; \tilde{z}_0= \left[\cos(x)\frac{\partial \tilde{\varphi}_M}{\partial x}\right]\bigg{|}_{t=0,\;x=x_0}, \\
 \tilde{\Gamma}_0 &= \bigg{\{}\cos(x) \frac{\partial}{\partial x}\Big{[}\cos(x) \frac{\partial \tilde{\varphi}_M}{\partial x} \Big{]}\bigg{\}}\bigg{|}_{t=0,\;x=x_0}, \\
 \;\;\; \tilde{A}_0&=\left\{{\cal \tilde{L}}\left[\cos(x) \frac{\partial \tilde{\varphi}_M}{\partial x} \right]\right\}\bigg{|}_{t=0,\;x=x_0},
 \end{split}
 \right.
 \end{eqnarray}
 where ${\cal \tilde{L}}$ is defined by (\ref{ellip:ope}).

Similarly, it is found that the optimal $c_0$ (i.e. the minimum $\cal E$) is $-1$.  In case of $c_0=-1$,  the homotopy approximations read
\begin{eqnarray}
\left\{
\begin{split}
\tilde{\varphi}_1(t,x)&=\sin(x)+t\cos(x);\\
\tilde{\varphi}_2(t,x)&=\sin(x)+t\cos(x)-\frac{1}{2}t^2\sin(x);\\
\tilde{\varphi}_3(t,x)&=\sin(x)+t\cos(x)-\frac{1}{2}t^2\sin(x)-\frac{1}{6}t^3\cos(x);\\
\cdots\quad&
\end{split}
\right.
\end{eqnarray}
It is found that the homotopy approximations $\tilde{\varphi}_M(t,x)$ are the Maclaurin series of the exact solution $u(t,x)=\sin(x+t)$ about $t$. Besides, as shown in Table \ref{case5}, the homotopy approximations at the initial moment $\tilde{y}_0$, $\tilde{z}_0$, $\tilde{\Gamma}_0$ and $\tilde{A}_0$ are exactly the same as the exact solutions $y_0$, $z_0$, $\Gamma_0$ and $A_0$, respectively. In addition, $\tilde{\cal E}$ quickly decreases to $7 \times 10^{-29}$ in just 0.28 seconds CPU times!   So, convergent results of this kind of 2FBSDEs can be rather efficiently obtained by means of the HAM, too.

\begin{table}[!htb]
\tabcolsep 0pt
\caption{The squared residual error, the used CPU time (second) and the relative error between the homotopy approximations $\tilde{y}_0$, $\tilde{z}_0$, $\tilde{\Gamma}_0$, $\tilde{A}_0$ and the exact solutions $ y_0$, $z_0$, $\Gamma_0$, $A_0$ of Eq.~(\ref{4:fbsde}), respectively. Here, $m$ denotes the order of approximation.}
\vspace*{-12pt}\label{case5}
\begin{center}
\def\temptablewidth{1\textwidth}
{\rule{\temptablewidth}{1pt}}
\begin{tabular*}{\temptablewidth}{@{\extracolsep{\fill}}ccccccc}
$m$   &$\tilde{\cal E}$ &  $\tilde{y}_0-y_{0}$ &$\tilde{z}_0-z_{0}$ &  $\tilde{\Gamma}_0-\Gamma_{0}$ & $\tilde{A}_0-A_{0}$   & time (s)  \\
\hline
4 &$1\times 10^{-4}$   &  0 & 0 &  0 &  0  &   0.05   \\
8  &$2\times 10^{-9}$ & 0 & 0 &  0 &  0 &   0.08    \\
12 &$3\times 10^{-15}$ &  0& 0 &  0&  0 &  0.13 \\
16  &$8\times 10^{-22}$ & 0 & 0 &  0 &  0  & 0.19 \\
20  &$7\times 10^{-29}$  & 0 & 0&  0&  0 &  0.28\\
\end{tabular*}
{\rule{\temptablewidth}{1pt}}
 \end{center}
 \end{table}

\subsection{An example of high-dimensional FBSDEs} %
Finally, let us consider the $d$-dimensional decoupled FBSDE
 \begin{equation}
 \frac{\partial u}{\partial t}+\sum_{i=1}^{d}b_i\frac{\partial u}{\partial x_i}+\frac{1}{2}\sum_{i=1}^{d} \sigma_{i,i}^2\cdot\frac{\partial^2 u}{\partial x_i^2 }+f=0,   \label{33:fbsde}
 \end{equation}
in which
\begin{equation}
\left\{
\begin{split}
&b_i=\frac{1}{d}\;x_i\;\mathrm{e}^{-x_i^2},\;\;\;\;\sigma_{i,i}=\frac{1}{d}\;\mathrm{e}^{-x_i^2},\;\;\;\;i=1,2,\cdots,d,\\
&f=-\sum_{i=1}^{d}x_i z_i+\frac{1}{d^2}y-\frac{1}{d^3}\sum_{i=1}^{d}\left[\left(x_i^2+\mathrm{e}^{-2x_i^2} \right)\prod_{\mbox{\tiny$\begin{array}{c}
k=1\\
k\neq i \end{array}$}}^{d}(x_k+t)\right]\\
&\quad\quad\quad\quad\quad\quad-\frac{1}{d}\sum_{i=1}^{d}\left[x_i^2\sum_{\mbox{\tiny$\begin{array}{c}
j=1\\
j\neq i \end{array}$}}^{d}\prod_{\mbox{\tiny$\begin{array}{c}
k=1\\
k\neq i\\
k\neq j\end{array}$}}^{d}(x_k+t)\right], \label{33:fbsde:boundary}
\end{split}
\right.
\end{equation}
with the exact solutions
\begin{equation}
\left\{
\begin{split}
&y_t=\frac{1}{d}\sum_{j=1}^{d}\left[x_{t,j}^2 \prod_{\mbox{\tiny$\begin{array}{c}
k=1\\
k\neq j\end{array}$}}^{d}(x_{t,k}+t)\right], \\
&z_{t,i}=\frac{1}{d^2}\mathrm{e}^{-x_{t,i}^2}\sum_{\mbox{\tiny$\begin{array}{c}
j=1\\
j\neq i\end{array}$}}^{d}\left[x^2_{t,j}\prod_{\mbox{\tiny$\begin{array}{c}
k=1 \\
k\neq i\\
k\neq j\end{array}$}}^{d}(x_{t,k}+t) \right]+\frac{2x_{t,i}}{d^2}\mathrm{e}^{-x_{t,i}^2}\prod_{\mbox{\tiny$\begin{array}{c}
k=1\\
k\neq i\end{array}$}}^{d}(x_{t,k}+t). \nonumber
\end{split}
\right.
\end{equation}
Without loss of generality, let us choose $T=1$, $x_{0, i}=1$, $i=1,2,\cdots,d$. Then the exact values at the initial moment are
\begin{equation}
y_0=1,\;\;\;\;\;\;\;\; z_{0,i}=\frac{d+1}{\mathrm{e}d^2},\;\;\;\;\;\;\;\;i=1,2,\cdots,d. \nonumber
\end{equation}
According to Eqs.~(\ref{1:general:fbsde})-(\ref{1:general:fbsde:boundary}), the corresponding PDE reads
\begin{eqnarray}
{\cal N}_6[u]&=&\frac{\partial u}{\partial t}+\frac{1}{2d^2}\sum_{i=1}^{d}\left(\mathrm{e}^{-2x_i^2}\cdot\frac{\partial^2 u}{\partial x_i^2}\right)+\frac{u}{d^2}-\frac{1}{d}\sum_{i=1}^{d}\left[x_i^2\sum_{\mbox{\tiny$\begin{array}{c}
j=1\\
j\neq i\end{array}$}}^{d}\prod_{\mbox{\tiny$\begin{array}{c}
k=1\\
k\neq i\\
k \neq j\end{array}$}}^{d}(x_k+t)\right]\nonumber \\
&& - \frac{1}{d^3}\sum_{i=1}^{d}\left[(x_i^2+\mathrm{e}^{-2x_i^2})\prod_{\mbox{\tiny$\begin{array}{c}
k=1\\
k\neq i\end{array}$}}^{d}(x_k+t)\right]=0,  \label{33:u}
\end{eqnarray}
subjects to the boundary condition
\begin{equation}
u(t,\textbf{x})\bigg{|}_{t=1}=\frac{1}{d}\sum_{j=1}^{d}\left[x_{j}^2 \prod_{\mbox{\tiny$\begin{array}{c}
k=1\\
k\neq j\end{array}$}}^{d}(x_{k}+1)\right],  \label{33:boundary}
\end{equation}
where 
\[     \textbf{x} = \Big{\{}  x_1, x_2, \cdots, x_d \Big{\}}   \]
is a vector, and ${\cal N}_6$ is a differential operator defined by (\ref{33:u}), respectively.       

In the frame of the HAM, we construct the following zeroth-order deformation equations
\begin{eqnarray}
(1-q)\frac{\partial}{\partial t}\Big{[}\Phi(t,\textbf{x},q) -\varphi_{0}(t,\textbf{x})\Big{]} =c_0 \; q \; {\cal N}_6 \Big{[}\Phi(t,\textbf{x},q)\Big{]}
,  \label{33:zeroth}
\end{eqnarray}
subject to the boundary condition
\begin{equation}
\Phi(t,\textbf{x},q)\bigg{|}_{t=1}=\frac{1}{d}\sum_{j=1}^{d}\left[x_{j}^2 \prod_{\mbox{\tiny$\begin{array}{c}
k=1\\
k\neq j\end{array}$}}^{d}(x_{k}+1)\right],   \label{33:zeroth:boundary}
\end{equation}
where $\varphi_0$ is the initial guess of $u$ that satisfies the boundary condition (\ref{33:boundary}).  
In order to satisfy the boundary condition (\ref{33:boundary}), we choose the initial guess
\begin{equation}
\varphi_0(t,\textbf{x})=\frac{1}{d}\sum_{j=1}^{d}\left[x_{j}^2 \prod_{\mbox{\tiny$\begin{array}{c}
k=1\\
k\neq j\end{array}$}}^{d}(x_{k}+1)\right].
\end{equation}
Similarly,  the above equations build a continuous variation $\Phi(t,\textbf{x},q)$ from the initial guess $\varphi_0$ at $q=0$ to the unknown solution $u$ at $q=1$.  Expand $\Phi(t,\textbf{x},q)$ into power series of $q$, i.e.
\begin{equation}
\Phi(t,\textbf{x},q)=\sum_{m=0}^{+\infty}\varphi_{m}(t,\textbf{x})\cdot q^m. \label{33:power:series}
\end{equation}
Substituting (\ref{33:power:series}) into the zeroth-order deformation equations (\ref{33:zeroth}) and (\ref{33:zeroth:boundary}), then equating the like-power of $q$, we have the high-order deformation equations
\begin{eqnarray}
\frac{\partial}{\partial t}\Big{[}\varphi_m-\chi_m\varphi_{m-1}\Big{]} = c_0 \; \delta_{6,m-1},
\end{eqnarray}
subject to the boundary condition
\begin{equation}
\varphi_m(t,\textbf{x})\bigg{|}_{t=1}=0,
\end{equation}
where
\begin{eqnarray}
\delta_{6,m-1}  &=& {\cal D}_{m-1} \bigg{\{} {\cal N}_6 \Big{[}\Phi(t,\textbf{x},q)\Big{]}  \bigg{\}}\nonumber\\
&=& \Bigg{\{} \frac{\partial \varphi_{m-1}}{\partial t}+\frac{\varphi_{m-1}}{d^2}  +\frac{1}{2d^2}\sum_{i=1}^{d}\left(\mathrm{e}^{-2x_i^2}\cdot\frac{\partial^2 \varphi_{m-1}}{\partial x_i^2}\right) \nonumber\\
&& -\frac{1}{d^3}\sum_{i=1}^{d}\Bigg{[}(x_i^2+\mathrm{e}^{-2x_i^2})\prod_{\mbox{\tiny$\begin{array}{c}
k=1\\ 
k\neq i\end{array}$}}^{d}(x_k+t)\Bigg{]} (1-\chi_m) \nonumber \\
&& -\frac{1}{d}\sum_{i=1}^{d}\Bigg{[}x_i^2\sum_{\mbox{\tiny$\begin{array}{c}
j=1\\
j\neq i\end{array}$}}^{d}\prod_{\mbox{\tiny$\begin{array}{c}
k=1\\
k\neq i\\
k \neq j\end{array}$}}^{d}(x_k+t)\Bigg{]} (1-\chi_m)\Bigg{\}}, 
\end{eqnarray}
and  $\chi_m$ and ${\cal D}_n$ are defined by (\ref{def:chi}) and (\ref{def:D}), respectively.   
Obviously,  it is easy to gain its solution
\begin{equation}
\varphi_m(t,\textbf{x}) = \chi_m \; \varphi_{m-1}(t,\textbf{x}) + c_0 \int^{t}_1 \delta_{6, m-1}(t,\textbf{x}) \mathrm{d} t, 
\end{equation}
step by step, starting from $m=1$, by means of computer algebra software Mathematica.

Similarly, at $q=1$, we have the $M$th-order homotopy-approximation of $u(t,\textbf{x})$:
\begin{equation}
\tilde{\varphi}_M(t,\textbf{x})=\sum_{m=0}^{M}\varphi_{m}(t,\textbf{x}), \label{33:homotopy:appro}
\end{equation}
which gives the approximations of the initial values:
\begin{equation}
\tilde{y}_0=\tilde{\varphi}_M(0,\textbf{x}_0),\;\;\;\tilde{z}_{0,i}=\bigg{[}\frac{1}{d}\;\mathrm{e}^{-x_i^2}\;\frac{\partial \tilde{\varphi}_M}{\partial x_i}\bigg{]}\bigg{|}_{t=0,\;\textbf{x}=\textbf{x}_{0}},\;\;\;i=1,2,\cdots,d. \nonumber
\end{equation}

\begin{figure}[t]
    \begin{center}
        \begin{tabular}{cc}
            \includegraphics[width=3in]{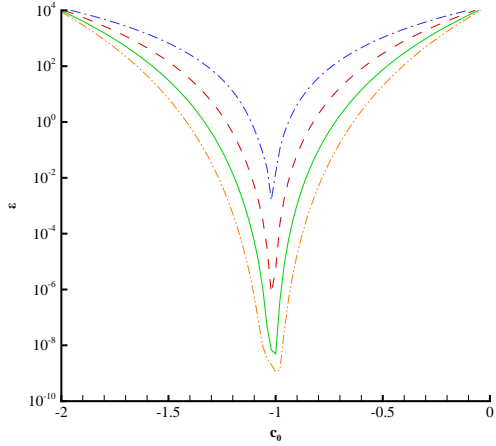}
        \end{tabular}
    \caption{The squared residual error  of Eq.~(\ref{33:u}) versus $c_0$  in case of $d=4$. Dash-dotted line: 2th-order approximation; Dashed line: 3th-order approximation; Solid line: 4th-order approximation; Dash-double-dotted line: 5th-order approximation.} \label{case6:c0}
    \end{center}
\end{figure}

\begin{figure}[!t]
    \begin{center}
        \begin{tabular}{cc}
            \includegraphics[width=3in]{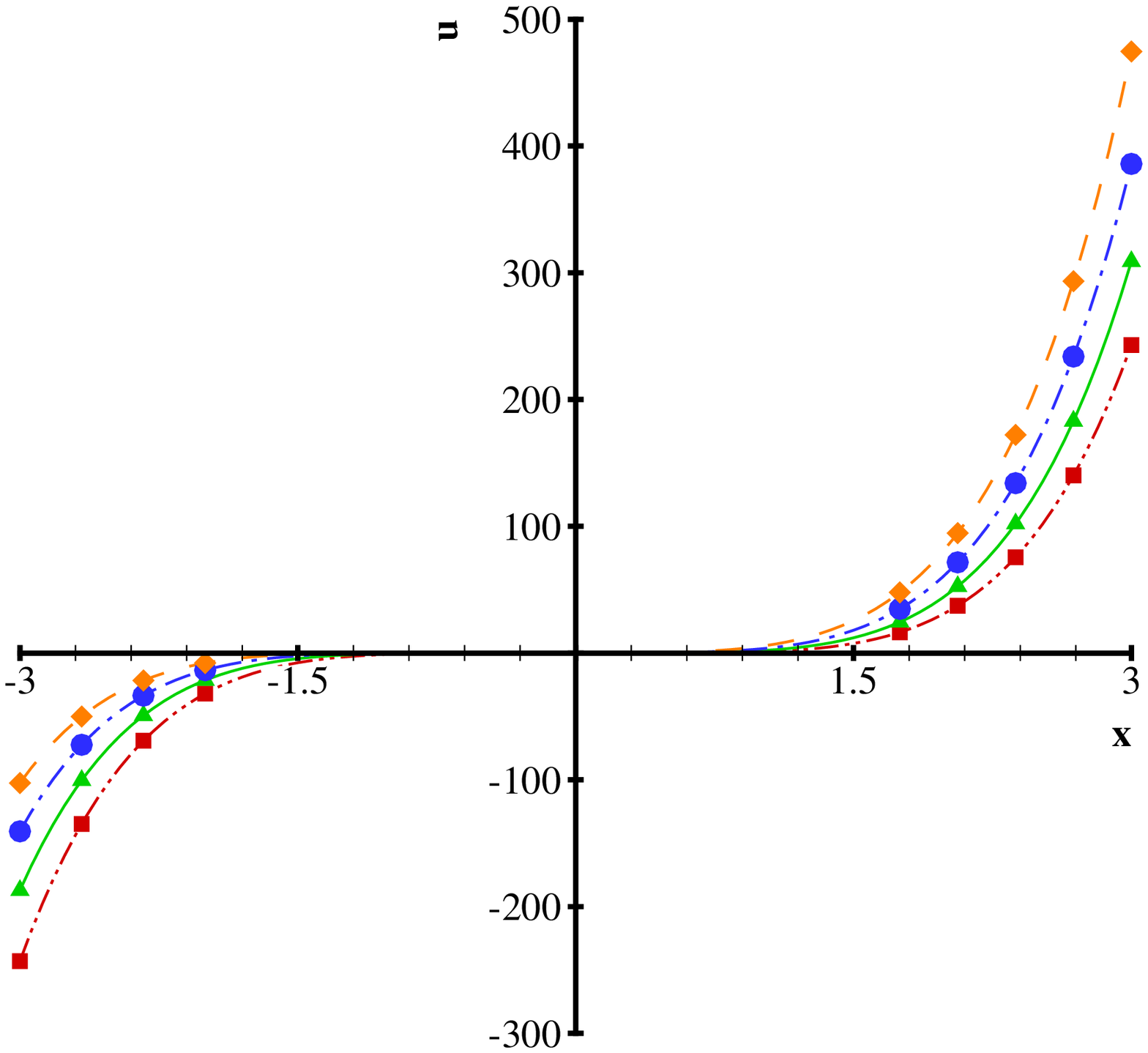}
        \end{tabular}
    \caption{Comparison between the homotopy approximations and the exact solution $u$ of Eqs.~(\ref{33:u}) and (\ref{33:boundary}) in case of $d=4$ when $t=0, 0.25, 0.5, 0.75$, $x_i=x$, $i=1,2,\cdots,d$. Dash-double-dotted line: homotopy approximation at $t=0$; Solid line: homotopy approximation at $t=0.25$; Dashed-dotted line: homotopy approximation at $t=0.5$; Dashed line: homotopy approximation at $t=0.75$. Square: exact solution at $t=0$; Triangle up: exact solution at $t=0.25$; Circle: exact solution at $t=0.5$; Diamond: exact solution at $t=0.75$.} \label{case6:u}
    \end{center}
\end{figure}

\begin{figure}[!h]
    \begin{center}
        \begin{tabular}{cc}
            \includegraphics[width=3in]{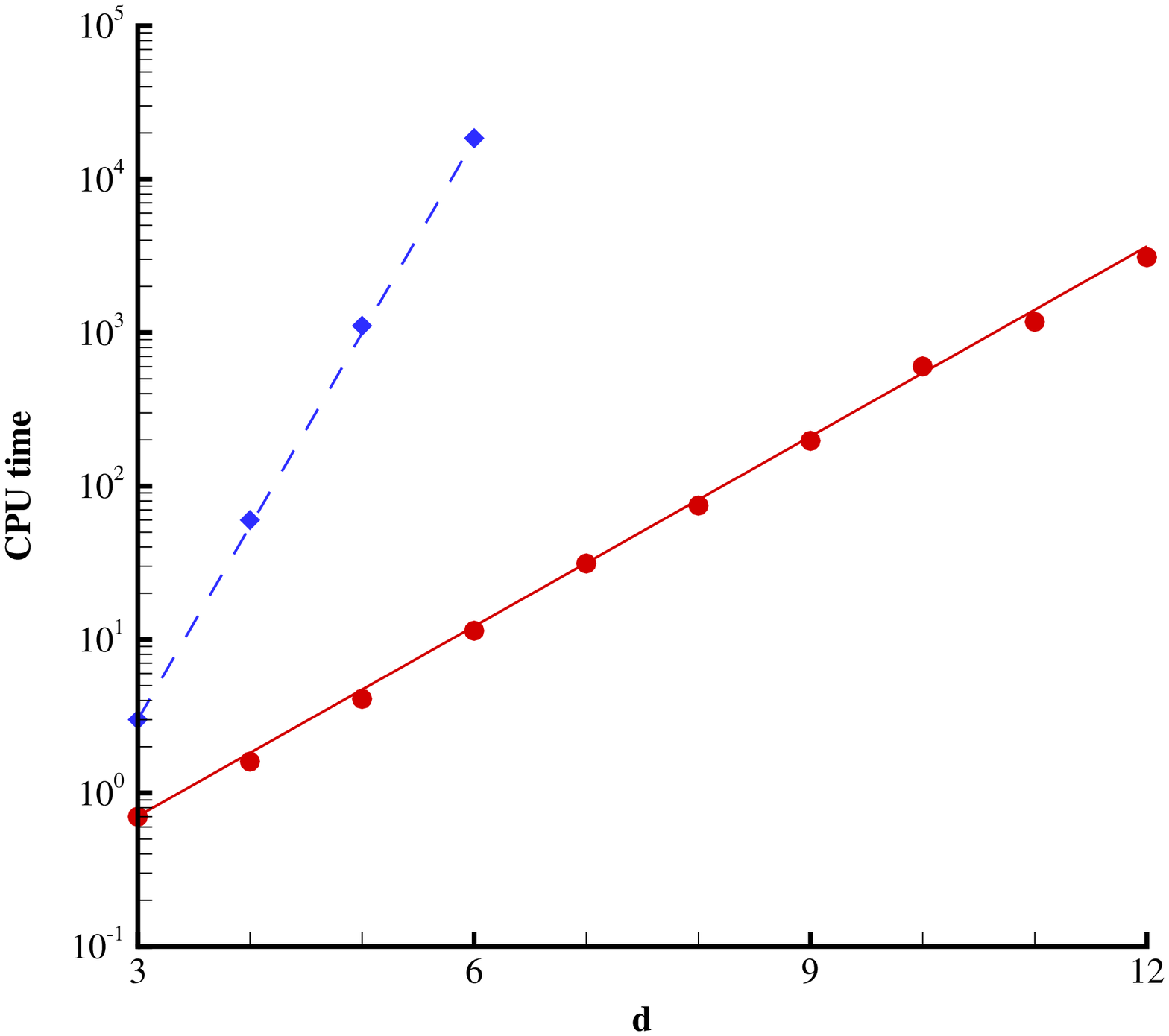}
        \end{tabular}
    \caption{The used CPU time (second) versus dimensionality $d$ of Eqs.~(\ref{33:u}) and (\ref{33:boundary}). Circle: CPU times for 5th-order HAM approximations; Diamond: CPU times for the spectral sparse grid approximations \cite{Fu2016Efficient};  Solid line: the fitted formula (\ref{def:t-H}); Dashed-line: the fitted formula (\ref{def:t-S}). } \label{CPU}
    \end{center}
\end{figure}

Similarly,  the curves of the residual squared error $\cal E$ versus $c_0$ in the case of $d=4$ is as shown in Fig.~\ref{case6:c0}, which indicates that the optimal $c_0$ (i.e. the minimum $\cal E$) is near to $-1$.  Besides, it is found that,  for \emph{any} $c_0\in[-1.1,-0.6]$, convergent results can be obtained.   So, for the sake of simplicity, let us choose $c_0=-1$.   As shown in Fig.~\ref{case6:u}, our homotopy approximations agree well with the exact solution.   It is found that  convergent results can be quickly obtained in a similar way, even though in case of high dimensionality $d$, such as $d=12$, as shown in Tables \ref{case6} to \ref{case13}.    Note that in case of the dimensionality $d=6$, it takes 18481 seconds CPU times (when $N=128$) with the accuracy at the level $10^{-7}$ by means of the so-called spectral sparse grid approximations \cite{Fu2016Efficient}.  However,  in the same case of the dimensionality $d=6$, it takes only 159 seconds CPU time by means of the HAM, which is less than $1\%$ of the CPU time used in \cite{Fu2016Efficient} for the same case, to obtain the 10th-order HAM approximation at the level of $10^{-13}$ that is much more accurate.   Note that, In case of the dimensionality $d=6$,  the HAM approach takes only 11 seconds to gain the 5th-order approximation at the level $10^{-8}$.   In fact, almost all of our 5th-order approximations given by the HAM reach the higher accuracy level than $10^{-7}$.   For example, even in the case of the dimensionality $d=12$,  our HAM approach takes only 3084 seconds to gain the 5th-order approximation at the accuracy level $10^{-9}$.   So, let $t_S$ and $t_H$ denote the used CPU times for  the spectral sparse grid approximations \cite{Fu2016Efficient} and the 5th-order HAM approximations, respectively.  We have the  fitted formulas:
\begin{eqnarray}
t_{H}&=&\exp(-3.2+0.95d),  \label{def:t-H} \\
t_{S}&=&\exp(-7.6+2.9d),   \label{def:t-S}
\end{eqnarray}
where the dimensionality $d\geq 3$ is an integer.   The computational efficiency versus the dimensionality $d$ is as shown in Fig.~\ref{CPU}.    Note that, as the dimensionality $d$ enlarges, the increase of computational complexity of our HAM approach is not as dramatic as the spectral sparse grid  method \cite{Fu2016Efficient}.  This  indicates  the validity and high efficiency of the HAM to solve FBSDEs with high dimensionality.

It should be emphasized that we solve the three kinds of FBDES simply via integrations with respect to the time $t$  in the frame of the HAM.   This is exactly the same as the BSDEs, as illustrated in \S~2.   Thus,  there are {\em no} fundamental differences in the frame of the HAM to solve either the BSDEs or the FBSDEs considered in this paper.       

\begin{table}[!htb]
\tabcolsep 0pt
\caption{The squared residual error, the used CPU time (second) and the relative error between the homotopy approximations $\tilde{y}_0$, $\tilde{z}_0$ and the exact solutions $ y_0$, $z_0$ of Eqs.~(\ref{33:fbsde}) and (\ref{33:fbsde:boundary}), respectively, when $d=4$. Here, $m$ denotes the order of approximation.}
\vspace*{-12pt}\label{case6}
\begin{center}
\def\temptablewidth{1\textwidth}
{\rule{\temptablewidth}{1pt}}
\begin{tabular*}{\temptablewidth}{@{\extracolsep{\fill}}ccccc}
$m$ &$\tilde{\cal E}$  &  $\tilde{y}_0-y_{0}$ &$\tilde{z}_0-z_{0}$    & time (s)  \\
\hline
2   &$1\times 10^{-3}$ &  $8\times 10^{-3}$ & $3\times 10^{-4}$    &   0.1   \\
4   &$1\times 10^{-10}$ & $2\times 10^{-6}$ & $4\times 10^{-7}$  &   0.7    \\
5   &$1\times 10^{-11}$ & $-2\times 10^{-7}$ & $8\times 10^{-8}$  &   1.6    \\
6  &$4\times 10^{-12}$ &  $-3\times 10^{-8}$& $3\times 10^{-9}$ & 3 \\
8  & $1\times 10^{-12}$ &$2\times 10^{-9}$ & $-7\times 10^{-10}$   &10 \\
10  & $9\times 10^{-13}$ &$-6\times 10^{-11}$ & $7\times 10^{-11}$ &   26\\
\end{tabular*}
{\rule{\temptablewidth}{1pt}}
 \end{center}
 \end{table}

 \begin{table}[!htb]
\tabcolsep 0pt
\caption{The squared residual error, the used CPU time (second) and the relative error between the homotopy approximations $\tilde{y}_0$, $\tilde{z}_0$ and the exact solutions $ y_0$, $z_0$ of Eqs.~(\ref{33:fbsde}) and (\ref{33:fbsde:boundary}), respectively, when $d=6$. Here, $m$ denotes the order of approximation.}
\vspace*{-12pt}\label{case8}
\begin{center}
\def\temptablewidth{1\textwidth}
{\rule{\temptablewidth}{1pt}}
\begin{tabular*}{\temptablewidth}{@{\extracolsep{\fill}}ccccc}
$m$  & $\tilde{\cal E}$&  $\tilde{y}_0-y_{0}$ &$\tilde{z}_0-z_{0}$    & time (s)  \\
\hline
2 &$2\times10^{-2}$  &  $9\times 10^{-3}$ & $2\times 10^{-4}$    &   0.6   \\
4 &$3\times10^{-11}$  & $5\times 10^{-7}$ & $5\times 10^{-8}$  &   5    \\
5 &$9\times10^{-15}$  & $-2\times 10^{-8}$ & $4\times 10^{-9}$  &   11    \\
6 &$3\times10^{-16} $&  $-1\times 10^{-9}$& $5\times 10^{-11}$ &  21 \\
8 &$5\times10^{-18} $& $1\times 10^{-11}$ & $-3\times 10^{-12}$   & 66 \\
10 &$2\times10^{-19}$ & $-1\times 10^{-13}$ & $6\times 10^{-14}$ &   159\\
\end{tabular*}
{\rule{\temptablewidth}{1pt}}
 \end{center}
 \end{table}

\begin{table}[!htb]
\tabcolsep 0pt
\caption{The squared residual error, the used CPU time (second) and the relative error between the homotopy approximations $\tilde{y}_0$, $\tilde{z}_0$ and the exact solutions $ y_0$, $z_0$ of Eqs.~(\ref{33:fbsde}) and (\ref{33:fbsde:boundary}), respectively, when $d=8$. Here, $m$ denotes the order of approximation.}
\vspace*{-12pt}\label{case10}
\begin{center}
\def\temptablewidth{1\textwidth}
{\rule{\temptablewidth}{1pt}}
\begin{tabular*}{\temptablewidth}{@{\extracolsep{\fill}}ccccc}
$m$  & $\tilde{\cal E}$& $\tilde{y}_0-y_{0}$ &$\tilde{z}_0-z_{0}$    & time (s)  \\
\hline
2 &1  &  $1\times 10^{-2}$ & $3\times 10^{-4}$    &   3   \\
4 &$1\times 10^{-10}$  & $3\times 10^{-7}$ & $2\times 10^{-8}$  &   34    \\
5 &$1\times 10^{-15}$  & $-6\times 10^{-9}$ & $5\times 10^{-10}$  &   74    \\
6 &$7\times 10^{-19}$ &  $-2\times 10^{-10}$& $3\times 10^{-12}$ &  149 \\
8 &$1\times 10^{-21}$  & $7\times 10^{-13}$ & $-8\times 10^{-14}$   & 451 \\
10 &$5\times 10^{-24}$  & $-2\times 10^{-15}$ & $5\times 10^{-16}$ &   1084\\
\end{tabular*}
{\rule{\temptablewidth}{1pt}}
 \end{center}
 \end{table}

 \begin{table}[!htb]
\tabcolsep 0pt
\caption{The squared residual error, the used CPU time (second) and the relative error between the homotopy approximations $\tilde{y}_0$, $\tilde{z}_0$ and the exact solutions $ y_0$, $z_0$ of Eqs.~(\ref{33:fbsde}) and (\ref{33:fbsde:boundary}), respectively, when $d=10$. Here, $m$ denotes the order of approximation.}
\vspace*{-12pt}\label{case12}
\begin{center}
\def\temptablewidth{1\textwidth}
{\rule{\temptablewidth}{1pt}}
\begin{tabular*}{\temptablewidth}{@{\extracolsep{\fill}}ccccc}
$m$ &$\tilde{\cal E}$  &  $\tilde{y}_0-y_{0}$ &$\tilde{z}_0-z_{0}$    & time (s)  \\
\hline
2 &58  &  $2 \times 10^{-2}$ & $4\times 10^{-4}$    &   21   \\
4  &$1 \times 10^{-9}$ & $2\times 10^{-7}$ & $8\times 10^{-9}$  &   196    \\
5  &$4 \times 10^{-15}$ & $-3\times 10^{-9}$ & $2\times 10^{-10}$  &   497    \\
6  &$2 \times 10^{-20}$&  $-6\times 10^{-11}$& $4\times 10^{-13}$ &  876 \\
8 &$2 \times 10^{-22}$ & $1\times 10^{-13}$ & $-6\times 10^{-15}$   &2920 \\
\end{tabular*}
{\rule{\temptablewidth}{1pt}}
 \end{center}
 \end{table}

  \begin{table}[!htb]
\tabcolsep 0pt
\caption{The squared residual error, the used CPU time (second) and the relative error between the homotopy approximations $\tilde{y}_0$, $\tilde{z}_0$ and the exact solutions $ y_0$, $z_0$ of Eqs.~(\ref{33:fbsde}) and (\ref{33:fbsde:boundary}), respectively, when $d=12$. Here, $m$ denotes the order of approximation.}
\vspace*{-12pt}\label{case13}
\begin{center}
\def\temptablewidth{1\textwidth}
{\rule{\temptablewidth}{1pt}}
\begin{tabular*}{\temptablewidth}{@{\extracolsep{\fill}}ccccc}
$m$ &$\tilde{\cal E}$  &  $\tilde{y}_0-y_{0}$ &$\tilde{z}_0-z_{0}$    & time (s)  \\
\hline
2 &4274  &  $5 \times 10^{-2}$ & $7\times 10^{-4}$    &   121   \\
3  &$1\times 10^{-2}$ & $1\times 10^{-4}$ & $2\times 10^{-6}$  &   424    \\
4  &$2\times 10^{-8}$&  $2\times 10^{-7}$& $6\times 10^{-9}$ &  1154 \\
5 &$3\times 10^{-14}$ & $-2\times 10^{-9}$ & $7\times 10^{-11}$   &3084 \\
\end{tabular*}
{\rule{\temptablewidth}{1pt}}
 \end{center}
 \end{table}

\section{Concluding remarks and discussions}  

In this paper, the homotopy analysis method (HAM) is successfully employed to solve the BSDEs, FBSDEs, 2FBSDEs and high dimensional FBSDEs (up to 12 dimension). By means of the HAM, which provides us great freedom to construct a continuous deformation and besides a simple way to guarantee the convergence of solution series, convergent results are quickly obtained for all equations considered in this paper. Besides, all of our results agree well with the exact solutions.   All of these demonstrate the validity and high efficiency of the HAM for the backward/forward-backward stochastic differential equations, especially in case of high dimensionality.   

Our HAM approach contains somethings unusual and different.  Note that, in the frame of the HAM,  all of the BSDE, FBSDE and 2FBSDE can be solved {\em without} any essential differences:   the governing equations for their high-order approximations are the same in essence, related with the auxiliary linear operator 
\begin{equation} 
 {\cal L}[u] = \frac{\partial u}{\partial t}. \label{def:L:HAM}
 \end{equation}
 In other words, all high-order homotopy approximations can be easily obtained by integration with respect to the time $t$.  This is very easy and quite efficient by means of computer algebra such as Mathematica.  This is the reason why our HAM approach is much more efficient than other numerical ones.  It should be emphasized that, it is the HAM that provides us great freedom to choose such a simple auxiliary linear operator, as mentioned by Liao  \cite{Liaobook, liaobook2}.   For example, using such kind of freedom of the HAM, one can solve the 2nd-order Gelfand differential equation even  by means of a $2d$-th order linear differential operator in the frame of the HAM, where $d=1,2,3$ is the dimensionality \cite{Liao2007-SAM}.  In addition, a nonlinear differential equation can be solved in the frame of the HAM even by means of directly defining an inverse mapping, i.e. without calculating any inverse operators at all \cite{Liao2016-DMMiM}.    It should be emphasized that, if perturbation method is used to solve Eqs. (\ref{21:tran:pde:gov}) and (\ref{21:tran:pde:boundary}), one had to handle governing equations with a more complicated linear operator
\begin{equation}
\hat{\cal L} [ u ]= \frac{\partial u}{\partial t} + \frac{\theta}{2}\left( 1-\theta \right)\left( 1 - 2 \theta\right) \frac{\partial u}{\partial \theta} 
+\frac{\theta^2}{2}\left( 1-\theta\right)^2 \frac{\partial^2 u}{\partial \theta^2}, \label{def:L:perturbation}
\end{equation}        
which are rather difficult to solve.   The above linear operator suggests that the time $t$ should be as important as other variables.   However, our HAM approach reveals that only the time $t$ is important for the BSDEs and FBSDEs, which has priority over others.     So,  the above linear operator (\ref{def:L:perturbation}) from perturbation techniques might mislead us greatly.   Fortunately, the linear operator (\ref{def:L:HAM}) used in our HAM approach  could tell us the truth.   

Without doubt, it is valuable to further apply the HAM to solve other types of BSDEs and FBSDEs in science, finance and engineering.
       
\section*{Acknowledgment}
Thanks to Prof. Shige Peng (Shandong University, China) for his introducing the high-dimensional FBSDEs to us.   This work is partly supported by National Natural Science Foundation of China (Approval No. 11272209 and 11432009).

\renewcommand{\theequation}{A.\arabic{equation}}
\setcounter{equation}{0}
\section*{Appendix}
Here, we use the mathematical induction to prove Eq.~(\ref{31:Maclaurin}), i.e.
\begin{equation}
\varphi_i=\frac{t^i}{i!}\sin(x+\frac{i}{2}\pi),\;\;\;\;i=0,1,2,\cdots.  \label{App:goal}
\end{equation}
When $i=0$, $\varphi_0=\sin(x)$ satisfies (\ref{App:goal}). Suppose when $0 \leq i \leq k$, we have
\begin{equation}
\varphi_i=\frac{t^i}{i!}\sin(x+\frac{i}{2}\pi)={\cal D}_i \left[\sin(tq+x)\right],\;\;\;\;i=0,1,2,\cdots, k, \label{App:known1}
\end{equation}
where ${\cal D}_i$ is defined by (\ref{def:D}). Besides, it is obvious that
\begin{eqnarray}
\begin{split}
{\cal D}_i \left[\cos(tq+x)\right]&~=~\frac{t^i}{i!}\cos\left(x+\frac{i \pi}{2} \right),\\
{\cal D}_i \left[\cos(2tq+2x)\right]&~=~\frac{(2t)^i}{i!}\cos\left(2x+\frac{i \pi}{2} \right). \label{App:known2}
\end{split}
\end{eqnarray}
In addition, according to Liao \cite{liaobook2}, we have
\begin{equation}
\begin{split}
&\quad~{\cal D}_{m}\left[\prod_{i=1}^{j}\alpha_{i}(q) \right]\\
&=\sum_{r_1=0}^{m}{\cal D}_{m-r_1}[\alpha_1] \sum_{r_2=0}^{r_1}{\cal D}_{r_1-r_2}[\alpha_2] \cdots\sum_{r_{j-1}=0}^{r_{j-2}}{\cal D}_{r_{j-2}-r_{j-1}}[\alpha_{j-1}]{\cal D}_{r_{j-1}}[\alpha_{j}].
\end{split}
\end{equation}
Substituting (\ref{App:known1}) and (\ref{App:known2}) into (\ref{31:highorder:delta}), we have
\begin{eqnarray}
\begin{split}
&\quad~\frac{\partial}{\partial t}\Big{(}\varphi_{k+1}-\chi_{k+1}\varphi_{k}\Big{)}\\
&=-\bigg{\{}\frac{\partial \varphi_{k}}{\partial t}-\frac{t^{k}}{k!}\cos\left(x+\frac{k}{2}\pi\right)\\
&\quad~+{\cal D}_{k}\Big{[}\sin(tq+x)\sin(tq+x)\cos(tq+x)\\
&\quad~+\frac{1}{2}\cos(tq+x)\sin(tq+x)\sin(tq+x)\sin(tq+x)\cos(tq+x)\\
&\quad~-\frac{1}{4}\cos(2tq+2x)\sin(tq+x)\sin(tq+x)\sin(tq+x)\\
&\quad~-\frac{1}{4}\sin(tq+x)\sin(tq+x)\sin(tq+x)\\
&\quad~-\cos(tq+x)\sin(tq+x)\sin(tq+x) \Big{]}
\bigg{\}}\\
&=-\bigg{\{}\frac{\partial \varphi_{k}}{\partial t}-\frac{t^{k}}{k!}\cos\left(x+\frac{k}{2}\pi\right)+{\cal D}_{k}\big{[}0\big{]}\bigg{\}}\\
&=\frac{t^{k}}{k!}\sin\left(x+\frac{k+1}{2}\pi \right)-\frac{\partial \varphi_{k}}{\partial t}.
\end{split}
\end{eqnarray}
Thus,
\begin{equation}
\frac{\partial \varphi_{k+1}}{\partial t}=(\chi_{k+1}-1)\frac{\partial \varphi_{k}}{\partial t}+\frac{t^{k}}{k!}\sin\left(x+\frac{k+1}{2}\pi \right)=\frac{t^{k}}{k!}\sin\left(x+\frac{k+1}{2}\pi \right). \nonumber
\end{equation}
Combining (\ref{31:highorder:boundary}), it is easy to obtain
\begin{equation}
\varphi_{k+1}=\frac{t^{k+1}}{(k+1)!}\sin\left(x+\frac{k+1}{2}\pi\right).
\end{equation}
So, (\ref{App:goal}) satisfies for any $k \in N$.
\section*{References}
\bibliographystyle{elsarticle-num}
\bibliography{finance}
\end{document}